\documentclass[11pt]{article}
\usepackage{amsfonts,mathrsfs}
\usepackage{multicol}
\usepackage{geometry}
\usepackage{float}
\usepackage{tikz,chemfig}
\usepackage{yhmath}
\usepackage{indentfirst}   
\usepackage{epic}
\usepackage{etex}
\usepackage{pinlabel}
\usepackage{ulem}
\usepackage{cases}              
\usepackage{setspace,color}
\usepackage{enumerate,mathdots}
\usepackage{amsmath,amssymb,bm} 
\usepackage{graphics}
\usepackage{epsf, epstopdf}
\usepackage{color, ifpdf}  
\usepackage{graphicx,booktabs,multirow}
\usepackage{latexsym}
\usepackage{appendix}

\allowdisplaybreaks
\numberwithin{equation}{section}

\usetikzlibrary{backgrounds}
\usetikzlibrary{arrows, patterns}
\usetikzlibrary{arrows.meta}
\tikzstyle{none}=[inner sep=0mm]
\tikzstyle{dotted}=[dash pattern=on \pgflinewidth off 2pt]
\tikzstyle{dashed}=[dash pattern=on 3pt off 3pt]
\usepackage{tkz-euclide}

\makeatletter
\pgfdeclarepatternformonly[\LineSpace]{my north east lines}{\pgfqpoint{-1pt}{-1pt}}{\pgfqpoint{\LineSpace}{\LineSpace}}{\pgfqpoint{\LineSpace}{\LineSpace}}%
{
	\pgfsetcolor{\tikz@pattern@color}
	\pgfsetlinewidth{0.4pt}
	\pgfpathmoveto{\pgfqpoint{0pt}{0pt}}
	\pgfpathlineto{\pgfqpoint{\LineSpace + 0.1pt}{\LineSpace + 0.1pt}}
	\pgfusepath{stroke}
}

\pgfdeclarepatternformonly[\LineSpace]{my north west lines}{\pgfqpoint{-1pt}{-1pt}}{\pgfqpoint{\LineSpace}{\LineSpace}}{\pgfqpoint{\LineSpace}{\LineSpace}}%
{
	\pgfsetcolor{\tikz@pattern@color}
	\pgfsetlinewidth{0.4pt}
	\pgfpathmoveto{\pgfqpoint{0pt}{\LineSpace}}
	\pgfpathlineto{\pgfqpoint{\LineSpace + 0.1pt}{-0.1pt}}
	\pgfusepath{stroke}
}
\makeatother

\newdimen\LineSpace
\tikzset{
	line space/.code={\LineSpace=#1},
	line space=8pt
}

\newcommand \tikzp[2]
{
	\begin{center}
		\begin{tikzpicture}[scale=#1]
			#2
		\end{tikzpicture}
	\end{center}
}

\newcommand \blue[1] {{\color{blue}#1}}

\newcounter{countcase}

\newcounter{countclaim}

\newcommand{\proof}{{\noindent {\em Proof}. }}

\newcommand{\proofend}{{\hfill$\Box$}
\setcounter{countclaim} {0}
\setcounter{countcase} {0}
}

\newcommand{\claimend}{{\hfill$\natural$}}

\newcommand \aln[2]
{
\begin{align}\label{#1}
#2
\end{align}
}

\newcommand \eqn[2]
{
	\begin{eqnarray}\label{#1}
		#2
	\end{eqnarray}
}

\newcommand{\N}{\mathbb{N}}
\newcommand{\cov}{\mathcal{H}}
\newcommand{\slop}{\mathcal{S}}
\newcommand{\SCRS}{\mathscr{S}}
\newcommand{\SCRE}{\mathscr{E}}

\newcommand{\SCRB}{\mathscr{B}}
\newcommand{\SCRNB}{\mathscr{\bar B}}
\newcommand{\SCRX}{\mathscr{X}}
\newtheorem{Cl}{Claim}

\newcommand{\C}{\mathcal{C}}
\newcommand{\g}{\ell}

\newcommand \Clm[2]
{
	\begin{Cl}\label{#1}
		#2
	\end{Cl}
}

\newtheorem{theo}{Theorem}[section]

\newtheorem{coro}[theo]{Corollary}

\begin{document}
		
	\title{DP color functions versus chromatic polynomials (II)}
	
	\author{Meiqiao Zhang\thanks{Email: nie21.zm@e.ntu.edu.sg and 
			meiqiaozhang95@163.com.}, 
		Fengming Dong\thanks{Corresponding author. Email: fengming.dong@nie.edu.sg 
		and donggraph@163.com.}
		\\
		\small National Institute of Education,
		Nanyang Technological University, Singapore
	}
	
	\date{}
	
	\maketitle
	\begin{abstract}
 For any connected graph $G$, let $P(G,m)$ and $P_{DP}(G,m)$ denote the chromatic polynomial and DP color function of $G$, respectively. 
 It is known that $P_{DP}(G,m)\le P(G,m)$ holds for every 
 positive integer $m$. 
 Let $DP_\approx$ (resp. $DP_<$) 
 be the set of graphs $G$ for which there exists an integer $M$ such that $P_{DP}(G,m)=P(G,m)$ 
 (resp. $P_{DP}(G,m)<P(G,m)$) 
 holds for all integers $m \ge M$.
Determining the sets $DP_\approx$ and $DP_<$ is a key problem on the study of the DP color function. 
For any edge set $E_0$ of $G$, let $\ell_G(E_0)$ be the
length of a shortest cycle $C$ in $G$ such that
$|E(C)\cap E_0|$ is odd whenever such a cycle exists, and $\ell_G(E_0)=\infty$ otherwise. Write $\ell_G(E_0)$ as $\ell_G(e)$ 
if $E_0=\{e\}$.

In this paper,  
we prove that if $G$ has a spanning tree $T$ such that 
$\ell_G(e)$ is odd for each $e\in E(G)\setminus E(T)$, 
the edges in $E(G)\setminus E(T)$ can be labeled as $e_1,e_2,\cdots, e_q$ with $\ell_G(e_i)\le \ell_G(e_{i+1})$ for all $1\le i\le q-1$ and each edge $e_i$ is contained in a cycle $C_i$ of length $\ell_G(e_i)$ 
with $E(C_i)\subseteq E(T)\cup \{e_j: 1\le j\le i\}$, then $G$ is a graph in $DP_{\approx}$.
As a direct application of this conclusion, all plane near-triangulations and complete multipartite graphs with at least three partite sets belong to $DP_{\approx}$.
We also show that if $E^*$ is an edge set of $G$ such that $\ell_{G}(E^*)$ is even and $E^*$ satisfies certain conditions, then $G$ belongs to $DP_<$.
In particular, if $\ell_G(E^*)=4$, where $E^*$ is a set of edges between two disjoint vertex subsets of $G$, then $G$ belongs to $DP_<$.
Both results extend known ones in [DP color functions versus chromatic polynomials, 
{\it Advances in Applied Mathematics} {\bf 134} (2022), article 102301].
	\end{abstract}

\section{Introduction}\label{secintro}

\subsection{Proper coloring, list coloring and DP coloring}
\label{subseccol}

In this article, we consider simple graphs only. 
For any graph $G$, let $V(G)$ and $E(G)$ be the vertex set and edge set of $G$, respectively. 
For any two disjoint subsets $V_1$ and $V_2$ of $V(G)$,
let $E_G(V_1,V_2)$ be the set of edges $uv\in E(G)$, 
where $u\in V_1$ and $v\in V_2$. 
For any non-empty subset $V_0$ of $V(G)$, let $G[V_0]$ denote the subgraph of $G$ induced by $V_0$. 
For any $A\subseteq E(G)$, 
let $V(A)$ be the set of vertices in $G$ which are incident to some edges in $A$, and let $G[A]$ be the subgraph of $G$ with vertex set $V(A)$ and edge set $A$ when $A\ne \emptyset$.
Let $G\langle A\rangle$ be the spanning subgraph of $G$ with edge set $A$
and $G-A=G\langle E(G)\setminus A\rangle$,
and denote by $c(A)$ the number of components of $G\langle A\rangle$. For any $u\in V(G)$, let $N_G(u)$ (or simply $N(u)$) be the set of the neighbors of $u$ in $G$.


Denote the set of positive integers by $\N$.
For any  $m\in \N$, 
let $[m]=\{1,\cdots,m\}$. 
For any graph $G$, a 
{\it proper coloring} of $G$ is a mapping $c:V(G)\rightarrow \N$, such that $c(u)\neq c(v)$ for all $uv\in E(G)$.
For any positive integer $m$, a 
{\it proper $m$-coloring} of $G$ is a proper coloring $c$ with $c(v)\in [m]$ for all $v\in V(G)$.
The \textit{chromatic polynomial} $P(G, m)$ of $G$ is a function which counts the number of proper $m$-colorings of $G$ for each $m \in \N$. 
The chromatic polynomial was originally designed as a tool to attack the Four Color Conjecture~\cite{Birk1912}, but later gained unique research significance because of its elegant properties, see~\cite{Dong2021, Dong2005, Read1988, Royle2009} for reference.

To generalize proper coloring, Vizing~\cite{Vizing1976} and Erd\H{o}s, Rubin and Taylor~\cite{Erdos1979} independently introduced the notion of list coloring. For any graph $G$, a \textit{list assignment} $L$ of $G$ is a mapping from $V(G)$ to the power set of $\N$, and an $L$-\textit{coloring} of $G$ is a proper coloring $c$ with $c(v)\in L(v)$ for all $v\in V(G)$. Denote the number of $L$-colorings of $G$ by $P(G, L)$.

$L$ is an \textit{$m$-list assignment} of $G$ if $|L(v)|=m$ holds for all $v\in V(G)$. Then the \textit{list color function} $P_l(G,m)$ of $G$ counts the minimum value of $P(G, L)$ among all $m$-list assignments $L$ for each $m\in\N$. Obviously, $P_l(G,m) \le P(G,m)$ holds for each $m\in \N$. And surprisingly, $P_l(G,m)=P(G,m)$ holds whenever $m> \frac{|E(G)|-1}{ln(1+\sqrt{2})}$ (see~\cite{Wang2017}).
While this implies that the list color function of some graph might not be a polynomial~\cite{Donner1992}, the list color function $P_l(G,m)$ of any graph $G$ inherits all the nice properties of its chromatic polynomial when $m$ is sufficiently large. See~\cite{Thomassen2009} for some open problems of list color functions.

To make breakthroughs in list coloring, Dvo\v{r}\'{a}k and Postle~\cite{Dvorak2018} recently defined the \textit{correspondence coloring}, or \textit{DP-coloring}. The formal definition is as follows.

For any graph $G$, a \textit{cover} of $G$ is an ordered pair $\cov=(L,H)$, where $H$ is a graph and $L$ is a mapping from $V(G)$ to the power set of $V(H)$ satisfying the conditions below:
\begin{itemize}
	\item the set $\{L(u):u\in V(G)\}$ is a partition of $V(H)$,
	\item for every $u\in V(G)$, $H[L(u)]$ is a complete graph,
	\item if $u$ and $v$ are not adjacent in $G$, then $E_H(L(u), L(v))=\emptyset$, and 
	\item for each edge $uv\in E(G)$, $E_H(L(u), L(v))$ is a matching.
\end{itemize}

For any cover $\cov=(L,H)$ of $G$,
$\cov$ is \textit{$m$-fold} 
if $|L(v)|=m$ for all $v\in V(G)$, and $\cov$ is \textit{full} if for each edge $uv\in E(G)$, $E_H(L(u), L(v))$ is a perfect matching.
An \textit{$\cov$-coloring} of $G$ is an independent set $I$
in $H$ with $|I|=|V(G)|$.
Obviously, any $\cov$-coloring $I$ of $G$ has the property that 
$|I\cap L(v)|=1$ for each $v\in V(G)$. 
Denote the number of $\cov$-colorings of $G$ by $P_{DP} (G,\cov)$.

The \textit{DP color function} $P_{DP}(G,m)$ of $G$, introduced by Kaul and Mudrock~\cite{KaulMudrock2021} in 2019,
 counts the minimum value of $P_{DP}(G,\cov)$ among all $m$-fold covers $\cov$ of $G$ for each $m\in\N$. Note that $P_{DP} (G,m) \le P_l(G,m)$ holds for each $m\in\N$. Therefore, for each $m\in\N$,
\aln{inPDPGM}
{
	P_{DP} (G,m) \le  P_l(G,m)\le P(G,m).
}
It is known that all the equalities in~(\ref{inPDPGM}) can hold simultaneously. For example, the authors of~\cite{KaulMudrock2021} proved that $P_{DP} (G,m) =P(G,m)$ holds for all $m\in \N$ when $G$ is a chordal graph.
However, different from list color functions, not the DP color functions of all graphs tend to be the same as their chromatic polynomials.
In~\cite{KaulMudrock2021}, it is shown that for any graph $G$ with even girth, there exists an $N \in\N$, such that $P_{DP}(G,m)<P(G,m)$ for all integers $m \ge N$. 
Therefore, how to characterize 
the two classes of graphs
$DP_\approx$ and $DP_<$
becomes a research focus in the study of DP color functions, where 
\begin{itemize}
	\item $DP_\approx$ is the set of graphs $G$ for which there exists an integer $M$ such that $P_{DP}(G,m)=P(G,m)$ holds for all integers $m \ge M$, and 
	\item $DP_<$ is the set of graphs $G$ for which there exists an integer $M$ such that $P_{DP}(G,m)<P(G,m)$ holds for all integers $m \ge M$. 
\end{itemize}

So far it is still unknown if there exists a graph $G$ such that 
$G\notin DP_\approx$ 
and $G\notin DP_<$.
Thus,  a characterization of the graphs in $DP_\approx$ or $DP_<$ does not necessarily guarantee a characterization of the graphs in the other class.

In this paper, we shall introduce 
our new findings on determining 
$DP_\approx$ and $DP_<$.

\subsection{Known results}
\label{subsecold}
Throughout this paper, we need only to consider connected graphs because for disconnected graph $G$ with components $G_1,\cdots,G_k$,
\aln{eprod}
{P_{DP}(G,m)=\prod _{i=1}^{k}P_{DP}(G_i,m).}

In this subsection, we introduce the known graphs contained in sets $DP_\approx$ and $DP_<$ respectively.

\newcommand{\covGm}{\cov_{G,m}}
\newcommand{\LGm}{L_{G,m}}
\newcommand{\HGm}{H_{G,m}}

Let $\covGm=(\LGm,\HGm)$
denote the special full $m$-fold cover 
such that $\LGm(u)=\{(u,i): i\in [m]\}$ for each vertex $u\in V(G)$ and $E_{\HGm}(\LGm(u),\LGm(v))=\{(u,i)(v,i): i\in [m]\}$ for each edge $uv\in E(G)$. 
Obviously, $P_{DP}(G,\covGm)=P(G,m)$ for all $m\in \N$.
Let $DP^*$ denote the set of graphs $G$ 
for which there exists $M\in \N$ such that for every $m$-fold cover $\cov=(L,H)$ of $G$, if $H\not\cong \HGm$, then $P_{DP}(G,\cov)> P(G,m)$ holds for all integers $m\ge M$. 
Apparently, $DP^*\subseteq DP_\approx$, but whether $DP^*= DP_\approx$ or not is currently unknown. 

On one hand, Mudrock and Thomason~\cite{MudrockThomason2021} 
showed that each graph with a dominating vertex belongs to $DP_{\approx}$.
Actually they proved that each graph with a dominating vertex belongs to $DP^*$.
Dong and Yang~\cite{Dong2022} then extended their conclusion to a large set of connected graphs
(see Theorem~\ref{Dong22th}).

Let $\C_G(e)$ be the set of cycles in $G$ containing $e$ with the minimum order. Note that $\C_G(e)=\emptyset$ if $e$ is a bridge. 
Denote by $\g_G(e)$ the \textit{girth of edge $e$ in $G$}, which is the order of any $C\in \C_G(e)$ if $\C_G(e)\neq \emptyset$; otherwise, $\g_G(e)=\infty$.

\begin{theo}[\cite{Dong2022}] \label{Dong22th}
	Let $G$ be a graph with a spanning tree $T$.
	 If  for each edge $e$ in $E(G)\setminus E(T)$, 
	 $\g_G(e)$ is odd and  
	 there exists 
	$C\in \C_G(e)$ such that 
	$\g_G(e') < \g_G(e)$ for each $e' \in E(C) \setminus (E(T) \cup \{e\})$, then $G\in DP^*$ and hence $G\in DP_\approx$.
\end{theo}

On the other hand, some families of graphs belonging to 
$DP_<$ were found. 
Kaul and Mudrock~\cite{KaulMudrock2021} discovered the fact that for any graph $G$ with an edge $e$,
if  $P(G-e, m) <m P(G, m)/(m -1)$, 
then $P_{DP}(G, m) <P(G, m)$ holds, and showed that 
every graph with an even girth belongs to $DP_<$.
The latter conclusion was extended to the following one.

\begin{theo}
	[\cite{Dong2022}]
	\label{Dong22DP<}
Graph $G$ belongs to $DP_<$ if $G$ contains an edge of even girth.
\end{theo}

An \textit{edge gluing} of vertex disjoint graphs $G_1$ and $G_2$ is a graph obtained by identifying an edge in $G_1$ and an edge in $G_2$ as a same one.
Then, it is easy to check~\cite{Dong2022,KaulMax21} that $G$ belongs to $DP_<$ if $G_1\in DP_<$ and either $G_1$ is a block of $G$ or $G$ is an edge-gluing of $G_1$ and some other graph. Therefore, as shown in~\cite{Dong2022}, Theorem~\ref{Dong22DP<} cannot be a characterization of all the graphs in $DP_<$ because by edge gluing any graph $G$ in $DP_<$ with a number of $3$-cycles, infinitely many graphs $G'$ in $DP_<$ can be obtained in which $\g_{G'}(e)=3$ holds for all $e\in E(G')$.

\subsection{New  results}\label{subsecnew}

In this article, we will further extend Theorems~\ref{Dong22th}
and~\ref{Dong22DP<}. We first give the definition of a family of graphs.

A graph $G$ is called \textit{DP-good} if $G$ has a spanning tree $T$ and a labeling $e_1,\cdots,e_q$ of the edges in $E(G)\setminus E(T)$,
where $q=|E(G)|-|E(T)|$, 
such that 
$\g_G(e_1)\le \cdots \le \g_G(e_q)$ and 
for each $i\in[q]$, $\g_G(e_i)$ is odd and $E(C_i)\subseteq E(T)\cup \{e_1,\cdots,e_{i}\}$ holds for some $C_i\in \C_G(e_i)$.
Obviously, the $q$ cycles $C_1,\cdots, C_q$ are pairwise distinct.

\begin{figure}[!ht]
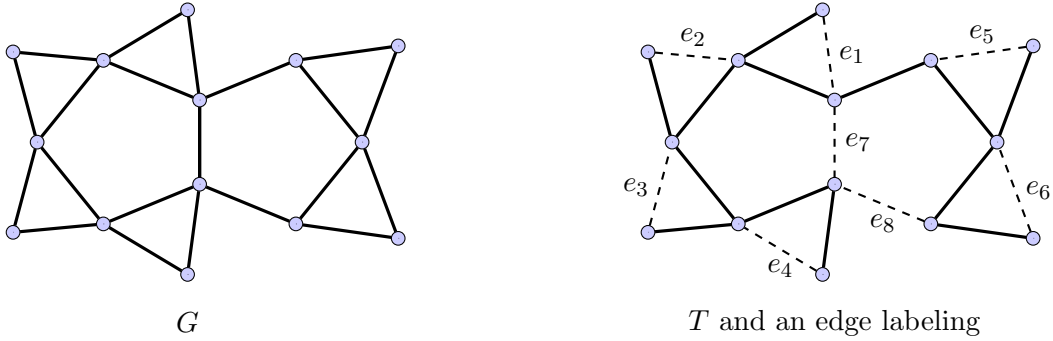

	\tikzstyle{cblue}=[circle, draw, thin,fill=blue!20, scale=0.5]
	\begin{multicols}{2}
		\tikzp{1.6}
		{
			\foreach \place/\y in {{(1.35,0)/1}, {(0.8,0.68)/2},{(0,0.35)/3},  {(0,-0.35)/4},{(0.8,-0.68)/5},{(-1.35,0)/6}, {(-0.8,0.68)/7},{(-0.8,-0.68)/8},
				{(-1.55, -0.75)/9},{(-1.55, 0.75)/10}, {(-0.1, -1.1)/11},{(-0.1, 1.1)/12},
				{(1.65, -0.8)/13},{(1.65, 0.8)/14}}   
			\node[cblue] (b\y) at \place {};
			
			\filldraw[black] (b1) circle (0pt)node[anchor=south] {};
			\filldraw[black] (b2) circle (0pt)node[anchor=east] {};
			\filldraw[black] (b3) circle (0pt)node[anchor=east] {};
			\filldraw[black] (b4) circle (0pt)node[anchor=north] {};
			\filldraw[black] (b5) circle (0pt)node[anchor=west] {};
			\filldraw[black] (b6) circle (0pt)node[anchor=west] {};
			\filldraw[black] (b7) circle (0pt)node[anchor=east] {};
			\filldraw[black] (b8) circle (0pt)node[anchor=west] {};		
			\filldraw[black] (b9) circle (0pt)node[anchor=west] {};
			\filldraw[black] (b10) circle (0pt)node[anchor=east] {};		
			\filldraw[black] (b11) circle (0pt)node[anchor=south] {};
			\filldraw[black] (b12) circle (0pt)node[anchor=south] {};
			\filldraw[black] (b13) circle (0pt)node[anchor=south] {};
			\filldraw[black] (b14) circle (0pt)node[anchor=south] {};
			
			\draw[black, very thick] (b1) -- (b2) -- (b3) -- (b4) -- (b5) -- (b1);
			\draw[black, very thick] (b3) -- (b7) -- (b6) -- (b8) -- (b4) -- (b3);
			\draw[black, very thick] (b6) -- (b9) -- (b8);
			\draw[black, very thick] (b6) -- (b10) -- (b7);
			\draw[black, very thick] (b8) -- (b11) -- (b4);
			\draw[black, very thick] (b7) -- (b12) -- (b3);
			\draw[black, very thick] (b5) -- (b13) -- (b1);
			\draw[black, very thick] (b1) -- (b14) -- (b2);
			
			\node [style=none] (cap1) at (-0.1, -1.5) {$G$};
		}
		\columnbreak
		\tikzp{1.6}
		{
			\foreach \place/\y in {{(1.35,0)/1}, {(0.8,0.68)/2},{(0,0.35)/3},  {(0,-0.35)/4},{(0.8,-0.68)/5},{(-1.35,0)/6}, {(-0.8,0.68)/7},{(-0.8,-0.68)/8},
				{(-1.55, -0.75)/9},{(-1.55, 0.75)/10}, {(-0.1, -1.1)/11},{(-0.1, 1.1)/12}, {(1.65, -0.8)/13},{(1.65, 0.8)/14}} 
			\node[cblue] (b\y) at \place {};
			
			\filldraw[black] (b1) circle (0pt)node[anchor=south] {};
			\filldraw[black] (b2) circle (0pt)node[anchor=east] {};
			\filldraw[black] (b3) circle (0pt)node[anchor=east] {};
			\filldraw[black] (b4) circle (0pt)node[anchor=north] {};
			\filldraw[black] (b5) circle (0pt)node[anchor=west] {};
			\filldraw[black] (b6) circle (0pt)node[anchor=west] {};
			\filldraw[black] (b7) circle (0pt)node[anchor=east] {};
			\filldraw[black] (b8) circle (0pt)node[anchor=west] {};		
			\filldraw[black] (b9) circle (0pt)node[anchor=west] {};
			\filldraw[black] (b10) circle (0pt)node[anchor=east] {};		
			\filldraw[black] (b11) circle (0pt)node[anchor=south] {};
			\filldraw[black] (b12) circle (0pt)node[anchor=south] {};
			\filldraw[black] (b13) circle (0pt)node[anchor=south] {};
			\filldraw[black] (b14) circle (0pt)node[anchor=south] {};
			
			\draw[black, thick, dashed]  (b9) -- node[left]{$e_3$}(b6);
			\draw[black, thick, dashed]  (b10) -- node[above]{$e_2$}(b7);
			\draw[black, thick, dashed]  (b11) -- node[below]{$e_4$}(b8);
			\draw[black, thick, dashed] (b12) -- node[right]{$e_1$}(b3);
			\draw[black, thick, dashed]  (b3)--node[right]{$e_7$}(b4);
			\draw[black, thick, dashed]  (b5)--node[below]{$e_8$}(b4);		
			\draw[black, thick, dashed]  (b14)--node[above]{$e_5$}(b2);		
			\draw[black, thick, dashed]  (b13)--node[right]{$e_6$}(b1);		
			\draw[black, very thick] (b5) -- (b1) -- (b2) -- (b3);
			\draw[black, very thick] (b3) -- (b7) -- (b6) -- (b8) -- (b4);
			\draw[black, very thick] (b8) -- (b9);
			\draw[black, very thick] (b6) -- (b10);
			\draw[black, very thick] (b4) -- (b11);
			\draw[black, very thick] (b7) -- (b12);
			\draw[black, very thick] (b5) -- (b13);
			\draw[black, very thick] (b1) -- (b14);
			
			\node [style=none] (cap1) at (0, -1.5) {$T$ and an edge labeling};
		}
	\end{multicols}
	\caption{A DP-good graph $G$ with a spanning tree $T$ and an edge labeling of the edges in $E(G)\setminus E(T)$}
	\label{goodfig}
\end{figure}

It is clear that any graph satisfying 
the condition in Theorem~\ref{Dong22th} is DP-good.
But the graph shown in Figure~\ref{goodfig} is a DP-good graph
which doesn't satisfy the requirement in Theorem~\ref{Dong22th}.
The following theorem shows 
that each DP-good graph belongs to $DP_\approx$.

\begin{theo}\label{theogd}
	Every DP-good graph is in $DP^*$.
\end{theo}

As an immediate consequence of Theorem~\ref{theogd}, Corollary~\ref{corogd4} below suggests that many special classes of graphs are DP-good and therefore contained in $DP^*$, such as chordal graphs, complete multipartite graphs with at least three partite sets, and plane near-triangulations.

\begin{coro}\label{corogd4}
	Let $G$ be a graph with vertex set $\{v_i: i=0,1,2,\cdots,n\}$.
	If for each $i\in [n]$, the set 
	$N(v_i)\cap \{v_j: 0\le j\le i-1\}$ is not empty and 
	the subgraph of $G$ induced by this vertex set is connected,
	then $G$ is DP-good.
\end{coro}

On the other hand, in order to extend Theorem~\ref{Dong22DP<}, we shall first generalize the definition of the girth of an edge to the girth of an edge set.
Given any subset $E_0$ of $E(G)$, let 
$\C'_G(E_0)$ be the set of the shortest cycles $C$ in $G$ such that $|E(C)\cap E_0|$ is odd (i.e., $|E(C)\cap E_0|$ is odd and $|E(C)\le |E(C')|$ holds for each cycle $C'$ in $G$ whenever $|E(C')\cap E_0|$ is odd), 
and {\it the girth of $E_0$}, 
denoted by $\g_G(E_0)$, is defined to be the length of any cycle in
$\C'_G(E_0)$ 
if this set is non-empty,  
and $\g_G(E_0)=\infty$ otherwise. 
Obviously, $\g_G(E_0)<\infty$
if and only if $G$ contains a cycle $C$ such that $|E(C)\cap E_0|$ is odd, 
and if $E_0=\{e\}$,  then $\C'_G(\{e\})=\C_G(e)$
and $\g_G(\{e\})=\g_G(e)$.

Let $E^*$ be a set of edges in $G$.
Assume that each edge $e$ in $E^*$ is assigned a direction $\overrightarrow{e}$
and $\overrightarrow{E^*}$ is the set of directed edges $\overrightarrow{e}$ for all $e\in E^*$. 
In graph $G$, only edges in $E^*$ are assigned directions. 
For any cycle $C$ in $G$,  
we say the directed edges 
in $\overrightarrow{E^*}$
are 
{\it balanced on $C$} if $|E(C)\cap E^*|$ is even and exactly half of the edges in $E(C)\cap E^*$ are 
oriented clockwise along $C$,
and \textit{unbalanced} otherwise. 
Obviously, the directed edges of 
$\overrightarrow{E^*}$ are balanced on $C$ 
when  $E(C)\cap E^*=\emptyset$,
and unbalanced on $C$ if 
$|E(C)\cap E^*|$ is odd.
Examples of cycles
on which directed edges of $\overrightarrow{E^*}$ are 
balanced or unbalanced
are shown in Figure~\ref{Ex10} (a) and (b), respectively, where $E(C)\cap E^*=\{e_1,e_2,e_3,e_4\}$.

\begin{figure}[!ht]
	\tikzstyle{cblue}=[circle, draw, thin,fill=blue!20, scale=0.5]
		\begin{multicols}{2}
		\tikzp{2.3}
		{
			\draw[very thick](0,0) ellipse (0.9 and 0.5);
			
			\foreach \place/\y in {{(-0.55,0.4)/1},{(-0.25,0.48)/2},{(0.25,0.48)/3},{(0.55,0.4)/4},{(-0.55,-0.4)/5},{(-0.25,-0.48)/6},{(0.25,-0.48)/7}, {(0.55,-0.4)/8},{(0.83,0.2)/9},{(0.83,-0.2)/10},{(-0.83,0.2)/11}, {(-0.83,-0.2)/12}}   
			\node[cblue] (b\y) at \place {};

			\filldraw[black] (b1) circle (0pt)node[anchor=south] {};
			\filldraw[black] (b2) circle (0pt)node[anchor=south] {};			
			\filldraw[black] (b3) circle (0pt)node[anchor=south] {};
			\filldraw[black] (b4) circle (0pt)node[anchor=south] {};
			\filldraw[black] (b5) circle (0pt)node[anchor=north] {};
			\filldraw[black] (b6) circle (0pt)node[anchor=north] {};
			\filldraw[black] (b7) circle (0pt)node[anchor=north] {};
			\filldraw[black] (b8) circle (0pt)node[anchor=north] {};
			
					\draw[black, -{>[scale=2.3,
	length=2, width=3]}, very thick] (-0.38,0.455) -- (-0.36,0.46);
\draw[black, -{>[scale=2.3,
	length=2, width=3]}, very thick] (0.38,0.455) -- (0.36,0.46);
					\draw[black, -{>[scale=2.3,
	length=2, width=3]}, very thick] (-0.43,-0.438) -- (-0.45,-0.432);
\draw[black, -{>[scale=2.3,
	length=2, width=3]}, very thick] (0.43,-0.438) -- (0.45,-0.432);
			
			\node [style=none] (cap1) at (1, 0.07) {$\vdots$};
			\node [style=none] (cap1) at (-1, 0.07) {$\vdots$};
			
			\node [style=none] (cap1) at (-0.43,0.6) {$\overrightarrow{e_1}$};
			\node [style=none] (cap1) at (0.41,0.6) {$\overrightarrow{e_2}$};
			\node [style=none] (cap1) at (-0.43,-0.63) {$\overrightarrow{e_4}$};
			\node [style=none] (cap1) at (0.41,-0.63) {$\overrightarrow{e_3}$};
		}
			\columnbreak
		\tikzp{2.3}
{
	\draw[very thick](0,0) ellipse (0.9 and 0.5);
	
			\foreach \place/\y in {{(-0.55,0.4)/1},{(-0.25,0.48)/2},{(0.25,0.48)/3},{(0.55,0.4)/4},{(-0.55,-0.4)/5},{(-0.25,-0.48)/6},{(0.25,-0.48)/7}, {(0.55,-0.4)/8},{(0.83,0.2)/9},{(0.83,-0.2)/10},{(-0.83,0.2)/11}, {(-0.83,-0.2)/12}}   
\node[cblue] (b\y) at \place {};

\filldraw[black] (b1) circle (0pt)node[anchor=south] {};
\filldraw[black] (b2) circle (0pt)node[anchor=south] {};			
\filldraw[black] (b3) circle (0pt)node[anchor=south] {};
\filldraw[black] (b4) circle (0pt)node[anchor=south] {};
\filldraw[black] (b5) circle (0pt)node[anchor=north] {};
\filldraw[black] (b6) circle (0pt)node[anchor=north] {};
\filldraw[black] (b7) circle (0pt)node[anchor=north] {};
\filldraw[black] (b8) circle (0pt)node[anchor=north] {};

\draw[black, -{>[scale=2.3,
	length=2, width=3]}, very thick] (-0.43,0.438) -- (-0.45,0.432);
\draw[black, -{>[scale=2.3,
	length=2, width=3]}, very thick] (0.43,0.438) -- (0.45,0.432);
\draw[black, -{>[scale=2.3,
	length=2, width=3]}, very thick] (-0.38,-0.455) -- (-0.36,-0.46);
\draw[black, -{>[scale=2.3,
	length=2, width=3]}, very thick] (0.43,-0.438) -- (0.45,-0.432);

\node [style=none] (cap1) at (1, 0.07) {$\vdots$};
\node [style=none] (cap1) at (-1, 0.07) {$\vdots$};

\node [style=none] (cap1) at (-0.43,0.6) {$\overrightarrow{e_1}$};
\node [style=none] (cap1) at (0.41,0.6) {$\overrightarrow{e_2}$};
\node [style=none] (cap1) at (-0.43,-0.63) {$\overrightarrow{e_4}$};
\node [style=none] (cap1) at (0.41,-0.63) {$\overrightarrow{e_3}$};
}
\end{multicols}

	{}\hfill \hspace{0cm} (a) 
	Balanced directed edges on $C$ \hspace{1.5 cm} (b) 
	Unbalanced directed edges on $C$\hfill {}
	
	\caption{
	 $E(C)\cap E^*=\{e_1,e_2,e_3,e_4\}$}
	\label{Ex10}
\end{figure}


We are now going to introduce the second main result in this article.

\begin{theo}\label{theoPDP<}
Let $G$ be a connected graph and $E^*$ be a set of edges in $G$.
If the following conditions are satisfied, then $G$ 
belongs to $DP_{<}$:
\begin{enumerate} 
\item $r_0=\g_G(E^*)$ is even; and 

\item
 there exists a way to assign an
 orientation $\overrightarrow{e}$ for  each edge $e\in E^*$ 
such that the directed edges in 
$\overrightarrow{E^*}=
\{\overrightarrow{e}:e\in E^*\}$ are balanced 
on each cycle $C$ of $G$ with 
$|E(C)|<r_0$.
\end{enumerate} 
\end{theo}


The following Corollary~\ref{coroPDP<}
of Theorem~\ref{theoPDP<} introduces a family of graphs in $DP_<$, including the graphs determined by Theorem~\ref{Dong22DP<}.

\begin{coro}\label{coroPDP<}
	Let $G$ be any graph and 
	let $E^*\subseteq E_G(V_1,V_2)$, where 
	$V_1$ and $V_2$ are disjoint vertex subsets of $V(G)$ 
	with $V_1\cup V_2\ne V(G)$.
	If the following conditions are satisfied, then $G\in DP_<$:
	\begin{enumerate}
		\item $r_0=\g_G(E^*)$ is even; and

		\item for each cycle $C$ in $G$
		such that $|E(C)\cap E^*|$ is positive, either 
		$|E(C)|\ge r_0$ or 
		no component of the subgraph $C- (E^*\cap E(C))$ is a $(v_1,v_2)$-path 
		for some $v_1\in V_1$ and $v_2\in V_2$.
		
	\end{enumerate}
\end{coro}

\begin{figure}[!ht]
	\tikzstyle{cblue}=[circle, draw, thin,fill=blue!20, scale=0.5]
	\tikzp{0.9}
	{
		\foreach \place/\y in {{(-0.65,-1)/1}, {(0.65,-1)/2},{(1.2,0)/3},  {(0.65,1)/4},
			{(-0.65,1)/5},{(-1.2,0)/6}, {(0,-1.8)/7},{(1.6,-0.8)/8},
			{(1.6,0.8)/9},{(0,1.8)/10}, {(-1.6,0.8)/11},{(-1.6,-0.8)/12}, {(1.4,-2.3)/13},{(2.9,0)/14},
			{(1.4,2.3)/15},{(-1.4,2.3)/16}, {(-2.9,0)/17},{(-1.4,-2.3)/18}} 
		\node[cblue] (b\y) at \place {};
		
		\filldraw[black] (b1) circle (0pt)node[anchor=south west] {\small{$u_1$}};
		\filldraw[black] (b2) circle (0pt)node[anchor=south east] {\small{$v_1$}};
		\filldraw[black] (b7) circle (0pt)node[anchor=north] {$u_2$};
		\filldraw[black] (b8) circle (0pt)node[anchor=north west] {\small{$v_2$}};		
		\filldraw[black] (b13) circle (0pt)node[anchor=north] {$v_3$};					
		\filldraw[black] (b18) circle (0pt)node[anchor=north] {$u_3$};	
		
		\draw[blue!70, very thick]  (b1) -- node[pos=0.55, below]{}(b2);
		\draw[blue!70, very thick]  (b7) --(b2) node[pos=0.2, below]{} ;
		\draw[blue!70, very thick]  (b7) -- node[pos=0.7, above]{}(b8);
		\draw[blue!70, very thick]  (b7) -- 
		node[pos=0.8, anchor=north]
		{}(b13);
		\draw[blue!70, very thick]  
		(b18) -- node[above]{}(b13);

		\draw[black, thick]  (b2) -- (b3) -- (b4) -- (b5) -- (b6) -- (b1); 
		\draw[black,  thick] (b2) -- (b8) -- (b3) -- (b9) -- (b4) -- (b10) -- (b5) -- (b11) -- (b6) -- (b12) -- (b1) -- (b7)  ; 	
		\draw[black, thick]  (b8) -- (b9) -- (b10) -- (b11) -- (b12) -- (b7); 
		\draw[black,  thick]  (b13) -- (b8) -- (b14) -- (b9) -- (b15) -- (b10) -- (b16) -- (b11) -- (b17) -- (b12) -- (b18) -- (b7) ; 			
		\draw[black, thick] (b13) -- (b14) -- (b15) -- (b16) -- (b17) -- (b18);

		\foreach \place/\y in {{(7.2,-0.2)/1}, {(8.4,-0.2)/2},{(6.6,-1.4)/3},  {(9,-1.4)/4},{(6.6,1.0)/5},{(9,1.0)/6}, {(5.6,-2.4)/7},{(10,-2.4)/8},
			{(5.6,2.4)/9},{(10,2.4)/10}} 
		\node[cblue] (z\y) at \place {};
		
		\filldraw[black] (z1) circle (0pt)node[anchor=east] {\small{$u_1$}};
		\filldraw[black] (z2) circle (0pt)node[anchor=west] {\small{$v_1$}};
		\filldraw[black] (z3) circle (0pt)node[anchor=east] {$u_2$};
		\filldraw[black] (z4) circle (0pt)node[anchor=west] {$v_2$};
		\filldraw[black] (z5) circle (0pt)node[anchor=east] {};
		\filldraw[black] (z6) circle (0pt)node[anchor=west] {};
		\filldraw[black] (z7) circle (0pt)node[anchor=east] {$u_3$};
		\filldraw[black] (z8) circle (0pt)node[anchor=west] {$v_3$};		
		\filldraw[black] (z9) circle (0pt)node[anchor=west] {};
		\filldraw[black] (z10) circle (0pt)node[anchor=east] {};

		\draw[blue!70, very thick]  (z3) -- node[below]{} (z4) ;
		\draw[blue!70, very thick]  (z3) -- (z8) node[pos=0.7, anchor=south]{};
		\draw[blue!70, very thick]  (z4) -- (z7) node[pos=0.7, anchor=south]{};
		\draw[blue!70, very thick]  (z4) -- (z1) node[pos=0.5, anchor=west]{};
		\draw[blue!70, very thick]  
		(z2) -- (z3) node[pos=0.5, anchor=east]{};

		\draw[black, thick] (z4) -- (z6) -- (z5) -- (z3); 	
		\draw[black, thick] (z3) -- (z1);
		\draw[black, thick] (z4) -- (z2);
		\draw[black, thick] (z1) -- (z5) -- (z7) -- (z9) -- (z5); 	
		\draw[black, thick] (z2) -- (z6) -- (z10) -- (z8) -- (z6); 	
		\draw[black, thick] (z3) -- (z7);
		\draw[black, thick] (z4) -- (z8);
		\draw[black, thick] (z9) -- (z6);
		\draw[black, thick] (z10) -- (z5);
		
	}
	
	{}\hfill \hspace{0.2 cm} (a) \hspace{6.2 cm} (b) \hfill {}
	
	\caption{Two graphs in $DP_<$}
	\label{Ex2}
\end{figure}

It is easy to verify that both the graphs in Figure~\ref{Ex2} satisfy the conditions in Corollary~\ref{coroPDP<} by taking $V_1=\{u_i:i=1,2,3\}$, $V_2=\{v_i:i=1,2,3\}$ and $E^*=E_G(V_1,V_2)$.
Note that the graph in Figure~\ref{Ex2} (b) belongs to a family of graphs stated in the following corollary,
which follows from 
Corollary~\ref{coroPDP<} directly.

\begin{coro}\label{coro2PDP<}
	Let $G$ be any graph and 
	let $E^*\subseteq E_G(V_1,V_2)$, where 
	$V_1$ and $V_2$ are disjoint vertex subsets of $V(G)$. 
	If  $\g_G(E^*)=4$, 
	then $G\in DP_<$.
\end{coro}

We will introduce some notations and 
fundamental results on an $m$-fold cover of a graph in Section~\ref{nota-mfold}.
We will then prove Theorem~\ref{theogd}
and Corollary~\ref{corogd4} 
in Section~\ref{secgood},
and Theorem~\ref{theoPDP<}
and Corollary~\ref{coroPDP<}
in Section~\ref{sec<}.
Finally, in Section~\ref{secplane}, 
we will apply Theorem~\ref{theoPDP<}
to determine 
some families of plane graphs 
belonging to $DP_<$.

\section{Notations and preliminary facts on an $m$-fold cover} 
\label{nota-mfold}

In this section, we introduce some 
notations and 
preliminary facts on an $m$-fold cover
which will be applied in the 
proofs of Theorems~\ref{theogd}
and~\ref{theoPDP<}.

Let $G$ be a graph.
By the definition of $P_{DP}(G,m)$, 
 $P_{DP}(G,m)$ is actually equal to the 
 minimum value of $P_{DP}(G,\cov)$'s 
 over all those full  
 $m$-fold covers $\cov=(L,H)$
 of $G$ with  
 $L(u)=\{(u,i): i=1,\cdots,m\}$ for every $u\in V(G)$.
 Now we assume that $\cov=(L,H)$ is any full $m$-fold cover of $G$ with $L(u)=\{(u,i): i=1,\cdots,m\}$ for every $u\in V(G)$.

 For any edge $e=uv$ in $E(G)$, let 
 $$X_e(G,\cov)=E_H(L(u),L(v))\setminus \{(u,i)(v,i): i=1,\cdots,m\}, \text{ and }$$
$$Y_e(G,\cov)=\{i\in[m]: (u,i)(v,j)\in X_e(G,\cov)\}.$$
Then $|X_e(G,\cov)|=|Y_e(G,\cov)|$, and if $(u,i)(v,j)\in X_e(G,\cov)$, $j$ is also included in set $Y_e(G,\cov)$ as $(u,j)(v,s)\in X_e(G,\cov)$ for some $s\neq j$.
 We say an edge $e$ in $G$ is \textit{horizontal} 
 with respect to $\cov$ if $X_e(G,\cov)=\emptyset$; \textit{sloping} otherwise.
 Denote the set of sloping edges in $G$ with respect to $\cov$ by $\slop_G(\cov)$.
For a given spanning tree $T$ of $G$, it is common to further assume that each edge in $T$ is horizontal with respect to $\cov$ because we can
rename the vertices in $L(u)$ for every vertex $u\in V(G)$ to guarantee that $E_H(L(u),L(v))=\{(u,i)(v,i): i=1,\cdots,m\}$ holds whenever $uv\in E(T)$, during which the structure of graph $H$ remains unchanged.

 \def \setg {{\cal G}}
 
Let $\SCRS(\cov)$ (simply  $\SCRS$) be the set of subsets $S$ of $V(H)$ with $|S\cap L(v)|=1$ 
 for each $v\in V(G)$.
 Clearly, $|S|=|V(G)|$ for each $S\in \SCRS$.
 For each $U\subseteq V(G)$, let $\SCRS|_U$ be the set of subsets $S$ of $V(H)$ such that $|S\cap L(v)|=1$ for each $v\in U$ and $S\cap L(v)=\emptyset$ for each $v\in V(G)\setminus U$. 
 Clearly, $|S|=|U|$ for each $S\in \SCRS|_U$, and 
  $\SCRS=\SCRS|_U$ when $U=V(G)$. 
  
For any subgraph $G_0$ of $G$, 
 let $H_{G_0}$ be the subgraph of $H$ 
 with vertex set $\cup_{u\in V(G_0)}L(u)$
 and edge set 
 $\cup_{uv\in E(G_0)}E_H(L(u),L(v))$.
 Let $\setg_{\cov}(G_0)$ be the set of graphs $H_{G_0}[S]$
 (i.e., the subgraph of $H_{G_0}$ induced by $S$),
 where  $S\in \SCRS|_{V(G_0)}$,
 such that $H_{G_0}[S]\cong G_0$.
 Note that $H_{G_0}[S]$ is the induced 
 subgraph $H[S]$ whenever $G_0$ is an induced subgraph of $G$.
 For each $j\in[m]$, let $S_j(G_0)=\{(v,j):v\in V(G_0)\}$ and write $H_{G_0}[S_j(G_0)]$ as $H_j[G_0]$.

For each edge $e=uv\in E(G)$, let $\SCRS_e$ be the set of $S\in \SCRS$ such that the two vertices in 
$S\cap (L(u)\cup L(v))$ are adjacent in $H$. 
For each $A\subseteq E(G)$, let 
$\SCRS_A=\cap_{e\in A}\SCRS_e.$ Then, by the inclusion-exclusion principle, 
\aln{eincluex}
{
	P_{DP}(G,\cov)=\sum_{A\subseteq E(G)}(-1)^{|A|}|\SCRS_A|,
}
which generalizes a well known property of the chromatic polynomial that
\aln{Pincluex}
{
	P(G,m)=\sum_{A\subseteq E(G)}(-1)^{|A|}m^{c(A)}.
}

For any graph $F$, 
let $\SCRB(F)$ be the set of bridges 
(i.e., cut-edges) in $F$, and 
let $\SCRNB(F)=E(F)\setminus \SCRB(F)$. 
Write $\SCRNB(G\langle A\rangle)$ as 
$\SCRNB(A)$ for any $A\subseteq E(G)$.
The following properties hold, as proved in~\cite{Dong2022}.
\begin{enumerate}
\item For any $A\subseteq E(G)$, if $G_1,G_2,\cdots G_{c(A)}$ are the components of $G\langle A\rangle$, then
\aln{epr-i}
{
	|\SCRS_A|=\prod\limits_{i=1}^{c(A)}|\setg_{\cov}(G_i)|.
}

\item For any connected subgraph $G_0$ of $G$, we have $|\setg_{\cov}(G_0)|\le m$, where the equality holds if
$\SCRNB(G_0)\cap \slop_G(\cov)=\emptyset$
(i.e., 
$\SCRNB(G_0)$ does not contain sloping edges with respect to $\cov$).

\item By Facts (i) and (ii), for each $A\subseteq E(G)$, we have $|\SCRS_A|\le m^{c(A)}$, where the equality holds if 
$\SCRNB(A)\cap \slop_G(\cov)=\emptyset$.

\item Let $\SCRE(\cov)$ (or simply $\SCRE$)  be the set of 
subsets $A$ of $E(G)$ such that $\SCRNB(A)$ contains at least one sloping edge with respect to $\cov$. Then Fact (iii) implies that 
\aln{epr-ii}
{
	P_{DP}(G,\cov)-P(G,m)=\sum_{A\in \SCRE}(-1)^{|A|}(|\SCRS_A|-m^{c(A)}).
}
\item Fact (iii) also implies that for any $k\in[n]$,
\aln{epr-iii}
{
	\sum_{\substack{A\in \SCRE\\ c(A)=k}}
(-1)^{|A|}(|\SCRS_A|-m^{c(A)})
\ge \sum_{\substack{A\in \SCRE,~c(A)=k\\ |A| \text{ is even}}}(|\SCRS_A|-m^{k})
}
and
\aln{epr-iv}
{
	\sum_{\substack{A\in \SCRE\\ c(A)=k}}
	(-1)^{|A|}(|\SCRS_A|-m^{c(A)})
	\le \sum_{\substack{A\in \SCRE,~c(A)=k\\ |A| \text{ is odd}}}(m^{k}-|\SCRS_A|).
}
\item For any $A\in \SCRE$ and any sloping edge $e$ in $ \SCRNB(A)$, let $G_1$ be the component of $G\langle A\rangle$ containing $e$.
Then $|V(G_1)|\ge \g_G(e)$ and $c(A)\le |V(G)|-\g_G(e)+1$, and
$|A|=\g_G(e)$ whenever $c(A)=|V(G)|-\g_G(e)+1$.
\end{enumerate}

\section{Proof of Theorem~\ref{theogd}.}\label{secgood}
Now we give the proof of Theorem~\ref{theogd}.

\noindent\textit{Proof of Theorem~\ref{theogd}.}
We need only to prove that there exists an $M\in\N$, such that whenever $m\ge M$, $P_{DP}(G,\cov)>P(G,m)$ holds for every full $m$-fold cover $\cov=(L,H)$ of $G$ with $H\not\cong \HGm$.

Suppose $n=|V(G)|$. As $G$ is DP-good, 
$G$ has a spanning tree $T$ and an edge labeling $e_1,\cdots,e_q$ of the edges in $E(G)\setminus E(T)$, such that $\g_G(e_1)\le \cdots \le \g_G(e_q)$ and for all $i\in[q]$, $\g_G(e_i)$ is odd and $E(C_i)\subseteq E(T)\cup \{e_1,\cdots,e_{i}\}$ for some $C_i\in \C_G(e_i)$.

Let $\cov=(L,H)$ be a full $m$-fold cover of $G$ with $L(u)=\{(u,i):i=1,\cdots,m\}$ for every $u\in V(G)$ and $H\not\cong \HGm$. 
We can further assume that $\slop_G(\cov)\neq \emptyset$ and all the edges in $E(T)$ are horizontal with respect to $\cov$. Then, every sloping edge $e$ in $G$ with respect to $\cov$ is of odd girth as $\slop_G(\cov)\subseteq E(G)\setminus E(T)$.

In the following, write $X_e(G,\cov)$ and $Y_e(G,\cov)$ simply as $X_e$ and $Y_e$ for any edge $e\in E(G)$.
Let $r=\min\{\g_G(e): e\in \slop_G(\cov)\}$, and let $E_0=\{e_{k_1},\cdots,e_{k_t}\}$ be the set of sloping edges in $G$ with $\g_G(e_{k_i})=r$, where $k_1<k_2<\cdots<k_t$.  Then $r$ is odd, $r\ge 3$ and $1\le t\le q$. 
Let $\SCRX_r=\cup _{i=1}^tX_{e_{k_i}}$. Then $|\SCRX_r|=\sum_{i=1}^t|X_{e_{k_i}}|\ge 1$.

Recall that the cycles $C_1,\cdots, C_q$ are pairwise distinct. We first prove the following three claims.

\Clm{cl1-2-1}
{$\sum\limits_{i=1}^t|\setg_{\cov}(C_{k_i})|\le mt-
	\left |\bigcup\limits_{i=1}^tY_{e_{k_i}}\right |$, i.e., 
$\sum\limits_{i=1}^t
(m-|\setg_{\cov}(C_{k_i})|)
\ge \left |\bigcup\limits_{i=1}^tY_{e_{k_i}}
\right |$.
}

\proof 
It suffices to prove the two facts below:
\begin{enumerate}
	\item $|\setg_{\cov}(C_{k_1})|=m-|Y_{e_{k_1}}|;$ and
	\item for any integer $p\in [t-1]$, $|\setg_{\cov}(C_{k_{p+1}})|\le m-|Y_{e_{k_{p+1}}}\setminus(\cup_{i=1}^pY_{e_{k_i}})|$.
\end{enumerate}

Since $E(C_{k_1})\subseteq E(T)\cup \{e_1,\cdots,e_{k_1}\}$ and $\g_G(e_1)\le\cdots\le\g_G(e_{k_1})=r$, $C_{k_1}$ contains exactly one sloping edge $e_{k_1}$.
Thus for any $j\in[m]$, $H_j[C_{k_1}-\{e_{k_1}\}]\cong  C_{k_1}-\{e_{k_1}\}$, and $H_j[C_{k_1}]\cong  C_{k_1}$ if and only if $j\notin Y_{e_{k_1}}$. Hence Fact (i) holds.

Similarly, for $p\in [t-1]$, all the edges in $E(C_{k_{p+1}})\setminus \{e_{k_1},e_{k_2},\cdots e_{k_{p+1}}\}$ are horizontal
as $E(C_{k_{p+1}})\subseteq E(T)\cup \{e_1,e_2,\cdots,e_{k_{p+1}}\}$ and $\g_G(e_1)\le\cdots\le\g_G(e_{k_{p+1}})=r$.
Let $j\in Y_{e_{k_{p+1}}}\setminus(\cup_{i=1}^pY_{e_{k_i}})$. 
Then, $H_j[C_{k_{p+1}}-\{e_{k_{p+1}}\}]\cong  C_{k_{p+1}}-\{e_{k_{p+1}}\}$ but $H_j[C_{k_{p+1}}]\not\cong  C_{k_{p+1}}$. Hence Fact (ii) holds and Claim~\ref{cl1-2-1} follows.
\claimend

\Clm{cl1-2-1'}
{The following inequality holds:
\aln{in2-0}
{
	\sum_{i=1}^t
	\left (m^{n-r+1}-|\SCRS_{E(C_{k_i})}|
	\right )
	\ge \frac{|\SCRX_r|}{q}m^{n-r}.
}
}
\proof 
Since $|Y_{e}|=|X_{e}|$ for every edge $e\in E(G)$, we have
\aln{in2-5}
{
	|\bigcup_{i=1}^tY_{e_{k_i}}|\ge \max_{i\in[t]}|X_{e_{k_i}}|\ge \frac{1}{t}\sum\limits_{i=1}^t|X_{e_{k_i}}|=\frac{1}{t}|\SCRX_r|\ge \frac{1}{q}|\SCRX_r|.
}
Then, by (\ref{epr-i}) and Claim~\ref{cl1-2-1},
\eqn{in2-3}
{
\sum_{i=1}^t(m^{n-r+1}-|\SCRS_{E(C_{k_i})}|)&=&\sum_{i=1}^t(m^{n-r+1}-m^{n-r}|\setg_{\cov}(C_{k_i})|)
	\nonumber\\
	&\ge&  |\bigcup_{i=1}^tY_{e_{k_i}}|m^{n-r}
		\nonumber\\
	&\ge & \frac{|\SCRX_r|}{q}m^{n-r}.
}
\claimend

\Clm{cl1-2-2} 
{The following inequality holds:
\aln{in2-1}
{
	\sum_{\substack{A\in \SCRE\\ c(A)=n-r+1}}(-1)^{|A|}(|\SCRS_A|-m^{c(A)})\ge \frac{|\SCRX_r|}{q}m^{n-r}.
}
}
\proof Recall that for any $A\in \SCRE$, $\SCRNB(A)$ contains a sloping edge $e$, where $\g_G(e)\ge r$. 
Thus, by (vi) in Section~\ref{nota-mfold}, $\g_G(e)=r=|A|$ holds whenever $c(A)=n-r+1$. Therefore,
\eqn{in2-2}
{
	\sum_{\substack{A\in \SCRE\\ c(A)=n-r+1}}(-1)^{|A|}(|\SCRS_A|-m^{c(A)})
		&=&\sum_{\substack{A\in \SCRE,~|A|=r\\ c(A)=n-r+1}}(-1)^{r}(|\SCRS_A|-m^{c(A)})
	\nonumber\\
	&=&\sum_{\substack{A\in \SCRE,~|A|=r\\ c(A)=n-r+1}}(m^{c(A)}-|\SCRS_A|),
}
where the last equality holds as $r$ is odd.

By (iii) in Section~\ref{nota-mfold}, $m^{c(A)}\ge |\SCRS_A|$ for any $A\subseteq E(G)$, hence
\eqn{in2-2'}
{
\sum_{\substack{A\in \SCRE\\ c(A)=n-r+1}}(-1)^{|A|}(|\SCRS_A|-m^{c(A)})
		 &\ge& \sum_{i=1}^t(m^{n-r+1}-|\SCRS_{E(C_{k_i})}|)
	\nonumber\\
	&\ge&  \frac{|\SCRX_r|}{q}m^{n-r},
}
where the last inequality follows from Claim~\ref{cl1-2-1'}.
\claimend

The rest of the proof is basically the same as in the proof of Theorem~\ref{Dong22th} that are given in~\cite{Dong2022}. 
For completeness, we restate  the proofs of Claims~\ref{Dongcl3}-\ref{Dongcl6} here with slight changes.

\Clm{Dongcl3}
{
	For any subgraph $G_0$ of $G$, if $\g_G(e)\le r$ for each sloping edge $e$ in $G_0$, then $|\setg_{\cov}(G_0)|\ge m-|\SCRX_r|$.
}
\proof
Since $\g_G(e)\le r$ for each sloping edge $e$ in $G_0$, 
each sloping edge in $G_0$ belongs to $E_0=\{e_{k_1},\cdots,e_{k_t}\}$.
Thus, for every $j\in [m]\setminus (\cup_{i=1}^tY_{e_{k_i}})$, $H_j[G_0]\cong  G_0$ holds, implying that
\aln{cle2-6}
{
	|\setg_{\cov}(G_0)|\ge m- \left |\bigcup_{i=1}^tY_{e_{k_i}}\right |\ge m- \sum_{i=1}^t|Y_{e_{k_i}}|=m- \sum_{i=1}^t|X_{e_{k_i}}|=m-|\SCRX_r|.
}
Hence  Claim~\ref{Dongcl3} holds.
\claimend

\Clm{Dongcl4}
{
	For any $A\in \SCRE$ with $c(A)=n-r$, we have $|\SCRS_A|\ge (m-|\SCRX_r|)m^{n-r-1}$.
}
\proof
Since $A\in \SCRE$, $\SCRNB(A)$ contains a sloping edge $e$ with $\g_G(e)\ge r$. Thus, by (vi) in Section~\ref{nota-mfold}, $G\langle A\rangle$ has a component $G_0$ with $e\in E(G_0)$ and $|V(G_0)|\ge r$. Moreover, as $c(A)=n-r$, $|V(G_0)|\le r+1$ holds, and for any other component $G'$ of $G\langle A\rangle$, $G'$ is either an isolated vertex or an edge, and thus $|\setg_{\cov}(G')|=m.$ 
Hence by (\ref{epr-i}), it suffices to prove that $|\setg_{\cov}(G_0)|\ge m-|\SCRX_r|.$

If $G_0$ is $2$-connected, then for every edge $e\in E(G_0)$, $\g_G(e)\le |V(G_0)|\le r+1$. Moreover, for each sloping edge $e$ in $G_0$, $\g_G(e)\le r$ as $r+1$ is even and $\g_G(e)$ is odd. Hence $|\setg_{\cov}(G_0)|\ge m-|\SCRX_r|$ holds by Claim~\ref{Dongcl3}. 

Otherwise, $G_0$ has exaxctly two blocks as $G_0$ contains a cycle $C$ with $|V(C)|\ge r$. Then, it is clear that the two blocks of $G_0$ are $G_0[V(C)]$ and an edge $f$, where $|V(C)|=r$ and $|\setg_{\cov}(G_0[f])|=m$. As $G_0[V(C)]$ is $2$-connected, for every edge $e\in E(G_0[V(C)])$, $\g_G(e)\le |V(C)|\le r$ holds, and thus $|\setg_{\cov}(G_0[V(C)])|\ge m-|\SCRX_r|$ follows from Claim~\ref{Dongcl3}. Consequently, $|\setg_{\cov}(G_0)|\ge m-|\SCRX_r|$ and
Claim~\ref{Dongcl4} holds.
\claimend

For any $s\in\N$ with $s\le n-r$, let $\phi_s$ be the number of subsets $A\subseteq E(G)$ such that $c(A)=s$, $G\langle A\rangle$ is not a forest and $|A|$ is even. 
\Clm{Dongcl5}
{The following inequality holds:
	\aln{cle2-7}
	{
		\sum_{\substack{A\in \SCRE\\ c(A)=n-r}}(-1)^{|A|}(|\SCRS_A|-m^{c(A)})\ge -\phi_{n-r}|\SCRX_r|m^{n-r-1}.
	}
}

\proof By (\ref{epr-iii}) and Claim~\ref{Dongcl4}, 
\eqn{cle2-8}
{
	\sum_{\substack{A\in \SCRE\\ c(A)=n-r}}
	(-1)^{|A|}(|\SCRS_A|-m^{c(A)})
	&\ge & \sum_{\substack{A\in \SCRE,~c(A)=n-r\\ |A| \text{ is even}}}(|\SCRS_A|-m^{n-r})
	\nonumber\\
	&\ge & \sum_{\substack{A\in \SCRE,~c(A)=n-r\\ |A| \text{ is even}}}(-|\SCRX_r|m^{n-r-1})
	\nonumber\\
	&\ge & -\phi_{n-r}|\SCRX_r|m^{n-r-1}.
}
\claimend

\Clm{Dongcl6}
{ For each $s\in [n-r-1]$, we have
	\aln{cle2-9}
	{
		\sum_{\substack{A\in \SCRE\\ c(A)=s}}(-1)^{|A|}(|\SCRS_A|-m^{c(A)})\ge -\phi_sm^{s}.
	}
}
\proof By (\ref{epr-iii}),
\eqn{cle210}
{
	\sum_{\substack{A\in \SCRE\\ c(A)=s}}
	(-1)^{|A|}(|\SCRS_A|-m^{c(A)})
	&\ge & \sum_{\substack{A\in \SCRE,~c(A)=s\\ |A| \text{ is even}}}(|\SCRS_A|-m^s)
	\nonumber\\
	&\ge & \sum_{\substack{A\in \SCRE,~c(A)=s\\ |A| \text{ is even}}}(-m^s)
	\nonumber\\
	&\ge & -\phi_sm^s.
}
\claimend

Now we are going to prove the main result by recalling (\ref{epr-ii}) that
$$P_{DP}(G,\cov)-P(G,m)=\sum_{A\in \SCRE}(-1)^{|A|}(|\SCRS_A|-m^{c(A)}).$$

By (vi) in Section~\ref{nota-mfold} and Claims~\ref{cl1-2-2},~\ref{Dongcl5},~\ref{Dongcl6}, we have
\eqn{cle211}
{
	P_{DP}(G,\cov)-P(G,m)&=&\sum_{s=1}^{n-r+1}\sum_{\substack{A\in \SCRE\\ c(A)=s}}(-1)^{|A|}(|\SCRS_A|-m^{c(A)})
	\nonumber\\
	&\ge&\frac{|\SCRX_r|}{q}m^{n-r}-\phi_{n-r}|\SCRX_r|m^{n-r-1}-\sum_{s=1}^{n-r-1}\phi_sm^s
	\nonumber\\
	&\ge &\frac{1}{q}m^{n-r}-\phi_{n-r}m^{n-r-1}-\sum_{s=1}^{n-r-1}\phi_sm^s,
}
where the last inequality holds when $m \ge q\phi_{n-r}$. As $q, \phi_1,\cdots,\phi_{n-r}$ are independent of the value of $m$, there exists $M_r\in\N$, such that $P_{DP}(G,\cov)-P(G,m)>0$ for all $m\ge M_r$. 
Let $M=\max\{M_r:3\le r\le n, r\text{ is odd}\}$. Then the result is proven.
\proofend

\vspace{0.5 cm}

The proof of Corollary~\ref{corogd4}
is given below.

\noindent\textit{Proof of Corollary~\ref{corogd4}.}
For $i=0,1,\cdots,n$, 
let $V_i=\{v_j: 0\le j\le i\}$ 
and $G_i=G[V_i]$.
Obviously, $G_n=G$.

By Theorem~\ref{theogd}, it suffices to show that for any $i\ge 1$,  
$G_i$ has a spanning tree $T_i$ 
and the edges in $E(G_i)\setminus E(T_i)$ can be labeled 
as $e_1,e_2,\cdots, e_{s_i}$ 
such that for all $j=1,2,\cdots, s_i$, 
$|V(C_j)|=3$ and 
$E(C_j)\subseteq 
E(T_i)\cup \{e_t:1\le t\le j\}$ hold 
for some $C_j\in \C_G(e_j)$.
 
The above conclusion is obvious for $i=1$, as $G_1\cong K_2$ by the given conditions.
Now assume that the above conclusion holds for $1\le i<n$.

Since $G[V_i\cap N(v_{i+1})]$ is connected, the vertices in 
$V_i\cap N(v_{i+1})$ can be labeled 
as $v_{q_0},v_{q_1},\cdots, v_{q_l}$,
where $l=|N(v_{i+1})\cap V_i|-1$,  
such that for any $1\le j\le l$, 
$N(v_{q_j})\cap \{v_{q_0},\cdots, v_{q_{j-1}}\}\ne \emptyset$.
Now, let $T_{i+1}$ be the spanning tree 
of $G_{i+1}$ obtained from $T_i$ by 
adding edge $v_{i+1}v_{q_0}$.
Let $s_{i+1}=s_i+l$, 
and for any $1\le j\le l$, 
let $e_{s_i+j}$ denote the edge 
$v_{i+1}v_{q_j}$. 
Then, it is obvious that 
for any  $1\le j\le l$,
$e_{s_i+j}$ is contained in a cycle 
$C_{s_i+j}$ of length $3$ 
such that $E(C_{s_i+j})\subseteq 
E(T_{i+1})\cup \{e_t:1\le t\le s_i+j\}$. 

 Hence the above conclusion holds for $i+1$. 
 Therefore $G_n$ is DP-good.
\proofend

By Corollary~\ref{corogd4}, chordal graphs, complete $k$-partite graphs, where $k\ge 3$, and plane near-triangulations are DP-good.

\section{Proof of Theorem~\ref{theoPDP<}
	\label{sec<}}
We shall prove Theorem~\ref{theoPDP<} in this section.

\noindent\textit{Proof of Theorem~\ref{theoPDP<}.}
Assume $|V(G)|=n$ and $E^*=\{e_1,\cdots,e_k\}$, where $k\ge 1$.

If $k=1$, then $\g_G(e_1)$ is even, and
the result follows from Theorem~\ref{Dong22DP<} directly.

In the following, we assume that $k\ge 2$.
For each $i\in [k]$, 
let $e_i$ be the edge $u_iv_i$ for $u_i,v_i\in V(G)$,
and let $\overrightarrow{e_i}$ be
the directed edge $(u_i,v_i)$  
with tail $u_i$.
By condition (ii) in Theorem~\ref{theoPDP<}, 
the directed edges in 
$\overrightarrow{E^*}
=\{\overrightarrow{e_i}: i\in [k]\}$
are balanced on 
every cycle $C$ in $G$ 
with $E(C)<r_0$.

For any positive integer $m$, let $\cov=(L, H)$ be the $m$-fold cover of $G$ defined below:
\begin{itemize}
	\item $L(x)=\{(x,i):i=1,\cdots,m\}$ for all $x\in V(G)$;
	\item $E_{H}(L(x),L(y))= \{(x,i)(y,i):i=1,\cdots,m\}$ for every edge $xy\in E(G)\setminus E^*$; and 
	\item $E_{H}(L(u_i),L(v_i))= \{(u_i,q)(v_i,q+1):q=1,\cdots,m-1\}\cup\{(u_i,m)(v_i,1)\}$ for every edge $e_i=u_iv_i\in E^*$.
\end{itemize}
Clearly, $\slop_G(\cov)=E^*$ (i.e., only edges in $E^*$ are sloping in $G$ with respect to $\cov$).

An {\it induced} cycle of $G$ is a cycle in $G$ which is induced by some subset of $V(G)$.
We first analyze the structure of  connected subgraphs $G_0$ of $G$ 
with $|V(G_0)|\le r_0$
by several claims.

\setcounter{Cl}{0}

\Clm{cl1-4-7'} 
{Let $C$ be a cycle in $G$.
	If $|V(C)|\le r_0$ and $|E(C)\cap E^*|$ is odd, then $C$ is an induced cycle of $G$ with $|V(C)|=r_0$.}

\proof  By Condition (i) in Theorem~\ref{theoPDP<},  $|V(C)|=r_0$ trivially holds.

Suppose that 
there exists $e\in E(G)\setminus E(C)$ such that $e$ joins two vertices in $V(C)$.
Then, $G$ contains a cycle $C'$ such that $V(C')\subseteq V(C)$,
$|V(C')|<|V(C)|=r_0$ and 
$|E(C')\cap E^*|$ is odd,
a contradiction to the definition of $r_0$.

Hence $G[V(C)]=C$
and Claim~\ref{cl1-4-7'} holds.
\claimend

\Clm{cl1-4-7} 
{Let $G_0=(V_0,E_0)$ be a $2$-connected subgraph of $G$.
	If $|V_0|\le r_0$ and $|E^*\cap E_0|=1$, then $|V_0|=r_0$ and $G_0$ is an induced cycle of $G$.}

\proof Since $G_0$ is $2$-connected,
$G_0$ contains a cycle $C$ with $|E(C)\cap E^*|=1$, where $|V(C)|\le |V_0|\le r_0$.
By Claim~\ref{cl1-4-7'}, $C$ is an induced cycle of $G$ with $|V(C)|=r_0$, which implies that $|V_0|=r_0$, $V_0=V(C)$ and $G_0$ is $C$.
Hence Claim~\ref{cl1-4-7} holds.
\claimend

\Clm{cl1-4-02}
{Let $G_0=(V_0,E_0)$ be a connected subgraph of $G$ with $|V_0|\le r_0$.
If  $|E^*\cap E_0|\ge 2$, then 
$G_0-(E^*\cap E_0)$ is disconnected.}

\proof 
Suppose that $G_0-(E^*\cap E_0)$ is connected. 
Let $e',e''\in E^*\cap E_0$,
and let $P$ be a path in 
$G_0-(E^*\cap E_0)$
connecting the two end-vertices of $e'$.
Consequently, the edge set $E(P)\cup\{e'\}$ forms a cycle $C$ in $G_0$ with $|V(C)|\le r_0$ and $E(C)\cap E^*=\{e'\}$.
By Claim~\ref{cl1-4-7'}, $C$ is an 
induced cycle with $|V(C)|=r_0$,
implying that $G_0$ is $C$,
a contradiction to the fact that $e''\in E_0$.
Hence Claim~\ref{cl1-4-02} holds.
\claimend

\Clm{cl1-4-02'}
{Let $G_0=(V_0,E_0)$ be a connected subgraph of $G$ with $|V_0|\le r_0$. If $|E^*\cap E_0|\ge 2$, then 
	no edge $e$ in $E^*\cap E_0$
	joins two vertices in any 
	 component $G'$ of $G_0-(E^*\cap E_0)$
 (i.e., each component $G'$ of $G_0-(E^*\cap E_0)$
 is an induced subgraph of $G_0$).}


\proof 
Assume that $G'$ is a component of $G_0-(E^*\cap E_0)$ and $e$ is an edge
in $E^*\cap E_0$ which joins two vertices in $G'$.

Then, $G'+e$ has a block, say $G_1$, which contains $e$. 
By Claim~\ref{cl1-4-02}, 
$|V(G')|<|V(G_0)|$.
Thus, $|V(G_1)|\le |V(G')|
<|V_0|\le r_0$.
But, as $|E(G_1)\cap E^*|=1$, 
Claim~\ref{cl1-4-7} implies that 
$|V(G_1)|=r_0$, 
a contradiction.
\claimend

\Clm{cl1-4-03}
{Let $G_0=(V_0,E_0)$ be a connected subgraph of $G$ with $|V_0|\le r_0$.
Assume that
$\{U_1,U_2\}$ is a  
partition of $V_0$
such that $E_0\cap E^*=E_{G_0}(U_1,U_2)$ and
$G_0[U_i]$ is connected for both $i=1,2$.
If $G_0$ is not a cycle 
of length $r_0$,  then for all the edges $e_{i_1},\cdots,e_{i_t}$ in $E_{G_0}(U_1,U_2)$, the vertices $u_{i_1},\cdots,u_{i_t}$ must be in the same set $U_s$ for some $s\in \{1,2\}$.
}

\begin{figure}[!ht]
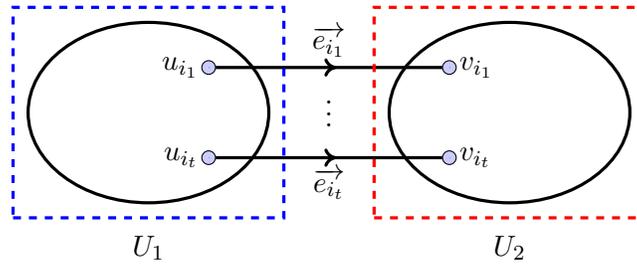

	\tikzstyle{cblue}=[circle, draw, thin,fill=blue!20, scale=0.5]
	\tikzp{2}
	{
		\draw[very thick](1.2,0) ellipse (0.8 and 0.6);
		\draw[very thick](-1.2,0) ellipse (0.8 and 0.6);

		\foreach \place/\y in {{(-0.8,0.3)/1}, {(0.8,0.3)/2},{(-0.8,-0.3)/3}, {(0.8,-0.3)/4}}   
		\node[cblue] (b\y) at \place {};
		
		\filldraw[black] (b1) circle (0pt)node[anchor=east] {$u_{i_1}$};
		\filldraw[black] (b2) circle (0pt)node[anchor=west] {$v_{i_1}$};
		\filldraw[black] (b3) circle (0pt)node[anchor=east] {$u_{i_t}$};
		\filldraw[black] (b4) circle (0pt)node[anchor=west] {$v_{i_t}$};		
		
		\draw[ very thick]  (b1) -- node[above]{$\overrightarrow{e_{i_1}}$}(b2);
		\draw[very thick]  (b3) -- node[below]{$\overrightarrow{e_{i_t}}$}(b4);
		
		\node [style=none] (vdots) at (0.0,0.05) {$\vdots$};
		
		\draw[black, -{>[scale=2.3,
			length=2, width=3]}, very thick] (0,0.3) -- (0.05,0.3);
		\draw[black, -{>[scale=2.3,
			length=2, width=3]}, very thick] (0,-0.3) -- (0.05,-0.3);
		
		\node [style=none] (cap1) at (-1.2,-0.9) {$U_1$};
		\node [style=none] (cap1) at (1.2,-0.9) {$U_2$};

		\draw[dashed, blue, very thick] (-2.1,-0.7) rectangle (-0.3,0.7);
		\draw[dashed, red, very thick] (2.1,-0.7) rectangle (0.3,0.7);
	}
	\caption{Graph $G_0$ with $E_{G_0}(U_1,U_2)=\{e_{i_1},\cdots, e_{i_t}\}$ and $u_{i_1},\cdots, u_{i_t}\in U_1$}
	\label{excl03}
\end{figure}

\proof 
Let $E'=E_{G_0}(U_1,U_2)$. If $|E'|=1$, then the result trivially holds.

In the following, assume that $|E'|\ge 2$. We need only to prove the two facts below on any two edges $e_{i_p},e_{i_q}$ in $E'$:
\begin{enumerate}
	\item if there is a cycle $C$ in $G_0$ shorter than $r_0$ with $|E(C)\cap E^*|= \{e_{i_p},e_{i_q}\}$, then $u_{i_p}$ and $u_{i_q}$ are contained in the same set $U_s$ for some $s\in \{1,2\}$;
	\item  
	otherwise, there exists $e_{i_j}\in E'\setminus \{e_{i_p},e_{i_q}\}$, such that there is a cycle $C_1$ in $G_0$ shorter than $r_0$ with $E(C_1)\cap E^*= \{e_{i_p},e_{i_j}\}$ and a cycle $C_2$ in $G_0$ shorter than $r_0$ with $E(C_2)\cap E^*= \{e_{i_j},e_{i_q}\}$.
\end{enumerate}

Since $G_0[U_i]$ is connected for both $i=1,2$, $G_0$ has a cycle $C$ 
with $|E(C)\cap E^*|= \{e_{i_p},e_{i_q}\}$.
If $|V(C)|<r_0$, Condition (ii) in Theorem~\ref{theoPDP<} indicates that 
$\overrightarrow{e_{i_p}}=(u_{i_p},v_{i_p})$ 
and $\overrightarrow{e_{i_q}}=(u_{i_q},v_{i_q})$ are balanced on $C$, 
implying that $u_{i_p}$ and $u_{i_q}$ must be in the same set $U_s$ for some $s\in \{1,2\}$. Fact (i) holds.

Now suppose that $G_0$ does not have
a cycle $C$ shorter than $r_0$ with $|E(C)\cap E^*|= \{e_{i_p},e_{i_q}\}$.
Thus,  $|V(C)|=r_0$, implying that 
$V(C)=V(G_0)$. 
As $G_0$ is 
not a cycle of length $r_0$, 
there is an edge $e\in E_0\setminus E(C)$.
Obviously, 
$e\notin E(G_0[U_1])\cup E(G_0[U_2])$.
Otherwise, $G_0$ has 
a cycle $C'$ shorter than $r_0$
with $|E(C')\cap E^*|= \{e_{i_p},e_{i_q}\}$,
a contradiction.
Thus,  $e\in E'=E_{G_0}(U_1,U_2)$.
Assume that $e=e_{i_j}\in E'\setminus \{e_{i_p},e_{i_q}\}$.
Then, there are cycles $C_1$ and $C_2$
in $C+e$  
such that $|E(C_1)\cap E^*|= \{e_{i_p},e_{i_j}\}$ and  $|E(C_2)\cap E^*|= \{e_{i_j},e_{i_q}\}$.
Note that both $C_1$ and $C_2$ are shorter than $r_0$.
Fact (ii) holds and Claim~\ref{cl1-4-03} follows.
\claimend

\Clm{cl1-4-04}
{Let $G_0=(V_0,E_0)$ be a connected subgraph of $G$ with $|V_0|\le r_0$ and 
$|E_0\cap E^*|\ge 2$. 
If $G_0$ is not a cycle of length $r_0$, then
 $|\setg_{\cov}(G_0)|=m$.}

\proof 
By Claim~\ref{cl1-4-02}, we can assume
that $G_0-(E_0\cap E^*)$ has $s$ 
($\ge 2$) components 
$G_1,G_2,\cdots,G_s$, 
where  $G_i=(V_i,E_i)$
for $i=1,2,\cdots,s$. 
Then, Claim~\ref{cl1-4-02'} implies that each $G_i$ is an induced subgraph of $G_0$, i.e., 
$E_0\cap E^*=
\bigcup_{1\le i<j\le s}
E_{G_0}(V_i,V_j)$.

Let $G'$ be the graph with vertex set  $V(G')=\{g_1,\cdots,g_s\}$
in which  $g_ig_j$ is an edge 
if and only if $E_{G_0}(V_i,V_j)\neq\emptyset$.
Let $\overrightarrow{G'}$ be the digraph obtained from $G'$ by 
converting each edge $g_ig_j$ in $G'$
into a directed edge 
whose tail is $g_i$ 
if and only if 
$u_q\in V_i$ for some edge 
$e_q=u_qv_q$ in $E_{G_0}(V_i,V_j)$. 
 Note that the orientation of directed edges in $\overrightarrow{G'}$ is well-defined due to the result in Claim~\ref{cl1-4-03}.
 An example is shown in Figure~\ref{figcl6-1} (b).

\begin{figure}[!ht]
	\tikzstyle{cblue}=[circle, draw, thin,fill=blue!20, scale=0.5]
	\tikzp{1.5}
	{
		\draw[fill=blue!6, very thick](-1.5,0) ellipse (0.55 and 0.4);
		\draw[fill=blue!6, very thick](0,0) ellipse (0.55 and 0.4);
		\draw[fill=blue!6, very thick](1.5,0) ellipse (0.55 and 0.4);
		\draw[fill=blue!6, very thick](-1.5,-2) ellipse (0.55 and 0.4);
\draw[fill=blue!6, very thick](0,-2) ellipse (0.55 and 0.4);
\draw[fill=blue!6, very thick](1.5,-2) ellipse (0.55 and 0.4);
		
		\foreach \place/\y in {{(4,0)/1}, {(5,0)/2},{(6,0)/3},{(4,-2)/4}, {(5,-2)/5},{(6,-2)/6}}   
		\node[cblue] (b\y) at \place {};
		
		\filldraw[black] (b1) circle (0pt)node[anchor=south] {$g_{1}$};
		\filldraw[black] (b2) circle (0pt)node[anchor=south] {$g_{2}$};
		\filldraw[black] (b3) circle (0pt)node[anchor=south] {$g_{3}$};
		\filldraw[black] (b4) circle (0pt)node[anchor=north] {$g_{4}$};		
		\filldraw[black] (b5) circle (0pt)node[anchor=north] {$g_{5}$};		
		\filldraw[black] (b6) circle (0pt)node[anchor=north] {$g_{6}$};				
		
		\draw[ very thick]  (b1) -- (b4) -- (b2) -- (b6) -- (b3) -- (b5) -- (b2);
		\draw[very thick]  (b3) -- (b4);
		
		\draw[black, -{>[scale=2.3, length=2, width=3]}, very thick] (4,-0.8) -- (4,-1);
		\draw[black, -{>[scale=2.3, length=2, width=3]}, very thick] (5,-0.5) -- (5,-0.7);
		\draw[black, -{>[scale=2.3, length=2, width=3]}, very thick] (6,-1) -- (6,-0.8);
		\draw[black, -{>[scale=2.3, length=2, width=3]}, very thick] (4.5,-1) -- (4.6,-0.8) ;
		\draw[black, -{>[scale=2.3, length=2, width=3]}, very thick] (4.5,-1.5) -- (4.6,-1.4);
		\draw[black, -{>[scale=2.3, length=2, width=3]}, very thick] (5.75,-0.5) -- (5.7,-0.6);
		\draw[black, -{>[scale=2.3, length=2, width=3]}, very thick] (5.75,-1.5) -- (5.7,-1.4);
				
		\draw[->, >=stealth, thick]  (-1.7,-0.25) -- (-1.7,-1.77);		
		\draw[->, >=stealth,thick]  (-1.6,-0.25) -- (-1.6,-1.76);
		\draw[->, >=stealth, thick]  (-1.5,-0.25) -- (-1.5,-1.75);

		\draw[->, >=stealth, thick]  (-1.25,-1.75) -- (1.3,-0.23);		
		\draw[ ->, >=stealth,thick]  (-1.35,-1.75) -- (-0.15,-0.25);
		\draw[->, >=stealth, thick]  (-1.45,-1.75) -- (-0.22,-0.2);
	
		\draw[->, >=stealth, thick]  (1.35,-1.75) -- (0.15,-0.25);
		\draw[->, >=stealth, thick]  (1.45,-1.75) -- (0.22,-0.21);
		\draw[->, >=stealth, thick]  (1.55,-1.75) -- (0.29,-0.17);		
						
		\draw[->, >=stealth, thick]  (-0.05,-0.25) -- (-0.05,-1.75);
		\draw[->, >=stealth, thick]  (0.05,-0.25) -- (0.05,-1.75);

		\draw[->, >=stealth, thick]  (1.7,-1.75) -- (1.7,-0.23);
		\draw[->, >=stealth, thick]  (1.6,-1.75) -- (1.6,-0.25);

		\draw[->, >=stealth, thick]  (1.5,-0.25) -- (0.22,-1.8);		
		\draw[->, >=stealth, thick]  (1.4,-0.25) -- (0.15,-1.75);

		\node [style=none] (cap1) at (-1.5,0)  {$G_1$};
		\node [style=none] (cap1) at (0,0)  {$G_2$};
		\node [style=none] (cap1) at (1.5,0)  {$G_3$};
		\node [style=none] (cap1) at (-1.5,-2)  {$G_4$};
		\node [style=none] (cap1) at (0,-2)  {$G_5$};
		\node [style=none] (cap1) at (1.5,-2)  {$G_6$};

	}
	{}\hfill \hspace{1.3 cm} (a) $G_0$ \hspace{6.3 cm} (b) $\overrightarrow{G'}$ \hfill {}

	\caption{An example of $\stackrel{\rightarrow}{G'}$}
	 \label{figcl6-1}
\end{figure}

As $G_0$ is connected, $G'$ is also connected.  
Let $T$ be a spanning tree of $G'$ with as many leaves as possible.
Thus, $T$ has at least $\Delta(G')$ leaves, where $\Delta(G')$ is the maximum degree of $G'$. 
For each vertex $g_i$ in $G'$, 
 there is a unique path, 
 denoted by $P_i$, in $T$ from $g_1$ to $g_i$. Denote by $\varphi_1(i)$ 
 the number of those edges 
 in $P_i$ whose
corresponding directed edges in 
$\overrightarrow{G'}$ are 
along the direction of path $P_i$ from $g_1$ to $g_i$, 
and denote by $\varphi_2(i)$ 
the number of the remaining edges in $P_i$.
Thus, $\varphi_1(i)+\varphi_2(i)=|E(P_i)|$.

\begin{figure}[!ht]
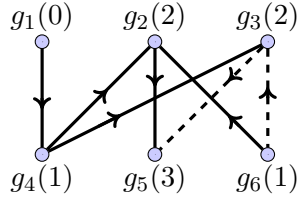

	\tikzstyle{cblue}=[circle, draw, thin,fill=blue!20, scale=0.5]
	\tikzp{1.5}
	{
		\foreach \place/\y in {{(4,0)/1}, {(5,0)/2},{(6,0)/3},{(4,-1)/4}, {(5,-1)/5},{(6,-1)/6}}   
		\node[cblue] (b\y) at \place {};
		
		\filldraw[black] (b1) circle (0pt)node[anchor=south] {$g_1 (0)$};
		\filldraw[black] (b2) circle (0pt)node[anchor=south] {$g_2 (2)$};
		\filldraw[black] (b3) circle (0pt)node[anchor=south] {$g_3 (2)$};
		\filldraw[black] (b4) circle (0pt)node[anchor=north] {$g_4 (1)$};		
		\filldraw[black] (b5) circle (0pt)node[anchor=north] {$g_5 (3)$};		
		\filldraw[black] (b6) circle (0pt)node[anchor=north] {$g_6 (1)$};				
		
		\draw[very thick]  (b1) -- (b4)-- (b2) -- (b5);		
		\draw[very thick]  (b4) -- (b3);
		\draw[very thick]  (b2) -- (b6);
		\draw[dashed, very thick]  (b3) -- (b6);
		\draw[dashed, very thick]  (b3) -- (b5);		
		\draw[black, -{>[scale=2.3, length=2, width=3]}, very thick] (4,-0.5) -- (4,-0.6);
		\draw[black, -{>[scale=2.3, length=2, width=3]}, very thick] (5,-0.3) -- (5,-0.4);
		\draw[black, -{>[scale=2.3, length=2, width=3]}, very thick] (6,-0.5) -- (6,-0.4);
		\draw[black, -{>[scale=2.3, length=2, width=3]}, very thick] (4.5,-0.75) -- (4.7,-0.65) ;
		\draw[black, -{>[scale=2.3, length=2, width=3]}, very thick] (4.5,-0.5) -- (4.58,-0.425);
		\draw[black, dashed, -{>[scale=2.3, length=2, width=3]}, very thick] (5.75,-0.75) -- (5.65,-0.65);
		\draw[black, dashed, -{>[scale=2.3, length=2, width=3]}, very thick] (5.7,-0.3) -- (5.65,-0.35);
	}
	\caption{With a spanning tree $T$
	consisting of dense edges, 
	the value of $\varphi(i)$ for each $i=1,2,\cdots,6$ is shown beside its veterx $g_i$}
	\label{figcl6-2}
\end{figure}

Now, let $\varphi(i)
=\varphi_1(i)-\varphi_2(i)$
for each $i\in [s]$. 
For the digraph $\overrightarrow{G'}$
in  Figure~\ref{figcl6-1} (b),
if $T$ is the spanning tree of $G'$ with its edge set $\{g_1g_4,g_2g_4,g_3g_4,g_2g_5,g_2g_6 \}$, then 
$$
\varphi(1)=0, \varphi(2)=\varphi(3)=2,
\varphi(4)=\varphi(6)=1,
\varphi(5)=3.
$$
as given in Figure~\ref{figcl6-2}.

We will complete the proof of this claim
by showing the following subclaims.

\noindent {\bf Subclaim~\ref{cl1-4-04}.1}. 
For any edge $g_ig_j\in E(T)$, 
$\varphi(j)=\varphi(i)+1$ whenever $(g_i,g_j)$ is the corresponding directed edge of $g_ig_j$ in 
$\overrightarrow{G'}$. 

Assume that $(g_i,g_j)$ is the corresponding directed edge of 
$g_ig_j$ in 
$\overrightarrow{G'}$.
As $g_ig_j\in E(T)$, either $g_i$ is on the path $P_j$, or $g_j$ is on the path $P_i$.
If $g_i$ is on the path $P_j$, then
$\varphi_1(j)=\varphi_1(i)+1$ and $
\varphi_2(j)=\varphi_2(i)$. 
If $g_j$ is on the path $P_i$, then
$\varphi_1(j)=\varphi_1(i)$ and 
$
\varphi_2(j)=\varphi_2(i)-1$. 
Thus, Subclaim~\ref{cl1-4-04}.1 follows
in both cases.

For every $q\in [m]$, 
let $S_q$ be the set in   
$\SCRS|_{V_0}$ defined as follows:
$$
S_q=\bigcup_{i=1}^s
\left \{(v,(q+\varphi(i))(\text{mod }m)):v\in V_i\right \},
$$
where $(v,0)=(v,m)$ for all $v\in V_0.$
Obviously, $\{S_1,\cdots,S_m\}$ is a partition of $V(H_{G_0})$.

\noindent {\bf Subclaim~\ref{cl1-4-04}.2}.
If $\varphi(j)=\varphi(i)+1$ holds 
for each directed edge 
$(g_i,g_j)$ in $\overrightarrow{G'}$,
then $H_{G_0}[S_q]\cong G_0$  for all $q\in [m]$, and hence Claim~\ref{cl1-4-04} holds.

Let $\phi$ be the bijection from $V_0$ to $S_q$ defined below:
for any 
$v\in V_0=\cup_{1\le i\le s}V_i$, 
$$
\phi(v)=(v,(q+\varphi(i))(\text{mod }m)),
\qquad \mbox{if }v\in V_i.
$$
To show that $H_{G_0}[S_q]\cong G_0$,
it suffices to prove that 
for each edge $uv\in E(G_0)$, 
$\phi(u)$ and $\phi(v)$ 
are adjacent in $H$.

For any $uv\in E_i$, 
where $1\le i\le s$,
we have $uv\in E_0\setminus E^*$, 
implying that $(u,(q+\varphi(i))(\text{mod }m))$
and 
$(v,(q+\varphi(i))(\text{mod }m))$
are adjacent in $H$ by the definition of $H$.
Now take any edge 
$uv\in E_{G_0}(V_i,V_j)\subseteq E^*$, 
where $1\le i,j\le s$.
Without loss of generality, 
assume that $(g_i,g_j)$ is an directed edge in 
$\overrightarrow{G'}$.
Then $\varphi(j)=\varphi(i)+1$ by the given condition in the subclaim.
By the definition of $H$,
$(u,(q+\varphi(i))(\text{mod }m))$
and 
$(v,(q+\varphi(j))(\text{mod }m))$
are adjacent in $H$.

Hence $H_{G_0}[S_q]\cong G_0$ for each $q\in [m]$, and the subclaim holds.

\noindent {\bf Subclaim~\ref{cl1-4-04}.3}.
$\varphi(j)=\varphi(i)+1$ holds 
for each directed edge 
$(g_i,g_j)$ in $\overrightarrow{G'}$.

Suppose that 
$\varphi(j)\ne \varphi(i)+1$ 
for some directed edge 
$(g_i,g_j)$ in $\overrightarrow{G'}$.
By Subclaim~\ref{cl1-4-04}.1,
$g_ig_j\in E(G')\setminus E(T)$.
Let $C'$ be the fundamental cycle of edge $g_{i}g_{j}$ in $G'$ 
with respect to spanning tree $T$. 
Assume that 
$g_{j_1},g_{j_2},\cdots,g_{j_t}$  
are the consecutive vertices on $C'$, where $t\ge 3$, $j_1=i$ and $j_t=j$.

As $G_i$ is connected for all $i=1,2,\cdots,s$, we can choose a 
shortest cycle $C$ in $G_0$ such that 
$$
E(C)\cap E^*\subseteq  \bigcup_{q=1}^t E_{G_0}(V_{j_{q}},V_{j_{q+1}}),
\quad 
|E(C)\cap E_{G_0}(V_{j_{q}},V_{j_{q+1}})|=1,
\quad \forall q\in [t],
$$ 
where $V_{j_{t+1}}=V_{j_1}$. 
Thus, $|E(C)\cap E^*|=t$.
Clearly, $t$ is an even integer; otherwise, 
Claim~\ref{cl1-4-7'} implies that 
$G_0$ is a cycle of length $r_0$, 
a contradiction. 

Suppose that $|V(C)|< r_0$.
By Condition (ii) in Theorem~\ref{theoPDP<},
the directed edges of 
$\overrightarrow{E^*}$ 
are  balanced
on $C$,
implying that 
the directed edges in  $\overrightarrow{G'}$ are balanced on $C'$.
By counting the number of edges in $C'$
which are oriented clockwise and counterclockwise along $C$ separately, 
we have 
$\varphi_1({j_1})+\varphi_2({j_t})+1=\varphi_1({j_t})+\varphi_2({j_1}),$
implying that  $\varphi(j_1)+1=\varphi(j_t)$
(i.e., $\varphi(i)+1=\varphi(j)$), 
a contradiction. 


Thus, $|V(C)|=r_0$, 
and so $V(C)=V_0$. 
Therefore, $t=s$ and $T$ is a path in $G'$ with $|V(T)|=|V(G')|=s$. 
Moreover, due to the choice of $C$, 
 for  each $q\in [s]$, 
$E_q\subseteq E(C)$
and $|E_{G_0}(V_{j_{q}},V_{j_{q+1}})|=1$,
implying that 
$E_{G_0}(V_{j_{q}},V_{j_{q+1}})
\subseteq E(C)$.

If $G'$ is a cycle,
then $G'$ is $C'$.
The above conclusion implies that 
$E_0=E(C)$, and thus 
$G_0$ is a cycle of length $r_0$, a contradiction.
Thus, $G'$ is not a cycle, implying that $G'$ has a spanning tree with at least three leaves,
a contradiction to the choice of $T$.

Hence Subclaim~\ref{cl1-4-04}.3 holds.

By Subclaims~\ref{cl1-4-04}.2 and~\ref{cl1-4-04}.3, 
$|\setg_{\cov}(G_0)|=m$ and Claim~\ref{cl1-4-04} holds.
\claimend

\Clm{cl1-4-9}
{For any $A\in \SCRE$,
if either $c(A)> n-r_0+1$
or  $c(A)=n-r_0+1$ and $|A|\neq r_0$,
then $|\SCRS_A|=m^{c(A)}$ holds.}

\proof As $A\in \SCRE$, $E^*\cap \SCRNB(A)\neq\emptyset$. Then, by (i) and (ii) in Section~\ref{nota-mfold}, it suffices to prove that for every block
$G_0=(V_0,E_0)$ of $G\langle A\rangle$ with $E_0\cap E^*\neq\emptyset$, $|\setg_{\cov}(G_0)|=m$ holds.

Suppose $G_0=(V_0,E_0)$ is a block of $G\langle A\rangle$ with $E_0\cap E^*\neq\emptyset$ and $|\setg_{\cov}(G_0)|< m$.
As $c(A)\ge n-r_0+1$, $|V_0|\le r_0$.
Then by Claims~\ref{cl1-4-7} and~\ref{cl1-4-04}, $|V_0|=|E_0|=r_0$, implying that 
either $c(A)< n-r_0+1$ or 
$c(A)= n-r_0+1$ and $|A|=r_0$, a contradiction. 
Hence Claim~\ref{cl1-4-9} holds.
\claimend

\Clm{cl1-4-20}
{If $m> k$,
then $\setg_\cov(C)=\emptyset$
for any cycle $C$ in $G$ 
such that $|E(C)\cap E^*|$ is odd.}

\proof Assume that $|E(C)\cap E^*|=2s+1$ for some integer $s\ge 0$ 
and $z_1,z_2,\cdots, z_q$
are consecutive vertices in $C$, where $q\ge 3$. 
Suppose that 
$\setg_\cov(C)\ne \emptyset$.
Then, there exists a cycle $C'$ in 
$H$ with consecutive vertices 
$(z_1,h_1), (z_2,h_2),\cdots,(z_q,h_q)$.
By the definition of $\cov=(L,H)$,
$h_{i+1}-h_i\ne 0$
if and only if $z_iz_{i+1}\in E^*$, 
and 
$
h_{i+1}-h_i\in \{0,1,-1,m-1,1-m\}
$
for all $i\in [q]$, 
where $h_{q+1}=h_1$ and $z_{q+1}=z_1$.
Thus, 
$h_{i+1}-h_i\ne 0$ holds for exactly  
$2s+1$ integers $i$'s in $[q]$.

Assume that there are exactly  $t$
integers   $i$'s in $[q]$ such that 
$h_{i+1}-h_i\in \{m-1,1-m\}$.
Then, there are exactly $(2s+1-t)$ 
integers   $i$'s in $[q]$ such that 
$h_{i+1}-h_i\in \{1,-1\}$.
It follows that 
\aln{eq1-cl8}
{
0=\sum_{i=1}^q (h_{i+1}-h_i)
=t'(m-1)+s'\times 1, 
}
where $t'$ and $s'$ are some integers 
with $|t'|\le t$ and 
$|s'|\le 2s+1-t$ 
such that 
both $t-t'$ and $(2s+1-t)-s'$
are even.

Suppose that $t'\ne 0$.
Without loss of generality, assume that 
$t'\ge 1$.
Then $s'\ge -(2s+1-1)=-2s$, and  
(\ref{eq1-cl8}) implies that 
\aln{eq2-cl8}
{
0=t'(m-1)+s'
\ge (m-1)-2s\ge (m-1)-(k-1)\ge m-k>0,
}
a contradiction. 
Hence $t'=0$, implying that $t$ is even. 
As 
$s'-(2s+1-t)$ is even, 
by (\ref{eq1-cl8}), 
\aln{eq3-cl8}
{
0=s'=(s'-(2s+1-t)) +(2s-t+1) 
\equiv 1 \pmod{2},
}
a contradiction.
Thus, Claim~\ref{cl1-4-20} holds.
\claimend

\Clm{cl1-4-10}
{The following inequality holds when $m> k$:
	\aln{in2-11}
	{
		\sum_{\substack{A\in \SCRE\\ c(A)=n-r_0+1}}
		(-1)^{|A|}(|\SCRS_A|-m^{c(A)})\le -m^{n-r_0+1}.}}
\proof
Let $C_0$ be any cycle in 
$\C'_G(E^*)$.
By Claim~\ref{cl1-4-20},  $|\SCRS_{E(C_0)}|=|\setg_\cov(C_0)|=0$
holds.

Obviously, $E(C_0)$ is a member in 
$\SCRE$ with $|E(C_0)|=r_0$ 
and $c(E(C_0))=n-r_0+1$.
Then, due to Claim~\ref{cl1-4-9}, the fact that $r_0$ is even, and (iii) in Section~\ref{nota-mfold}, we have 
\eqn{in2-12}
{
	\sum_{\substack{A\in \SCRE\\ c(A)=n-r_0+1}}
	(-1)^{|A|}(|\SCRS_A|-m^{c(A)})
	&=& \sum_{\substack{A\in \SCRE,~|A|=r_0\\ c(A)=n-r_0+1}}
	(-1)^{r_0}(|\SCRS_A|-m^{c(A)})
	\nonumber\\
	&=& \sum_{\substack{A\in \SCRE,~|A|=r_0\\ c(A)=n-r_0+1}}
	(|\SCRS_A|-m^{c(A)})
	\nonumber\\
	&\le& |\SCRS_{E(C_0)}|-m^{n-r_0+1}
	\nonumber\\
	&=& -m^{n-r_0+1}.
}
\claimend

For any $s\in\N$ with $s\le n-r_0$, let $\phi_s$ be the number of subsets $A\subseteq E(G)$ such that $c(A)=s$, $G\langle A\rangle$ is not a forest and $|A|$ is odd. 

\Clm{cl1-4-11} 
{For each $s\in [n-r_0]$, the following inequality holds:
	\aln{in2-15}
	{
		\sum_{\substack{A\in \SCRE\\ c(A)=s}}
		(-1)^{|A|}(|\SCRS_A|-m^{c(A)})\le \phi_sm^s.}
}
\proof By (\ref{epr-iv}),
\eqn{in2-16}
{
	 \sum_{\substack{A\in \SCRE\\ c(A)=s}}
	(-1)^{|A|}(|\SCRS_A|-m^{c(A)})
	&\le & \sum_{\substack{A\in \SCRE,~c(A)=s\\ |A| \text{ is odd}}}(m^s-|\SCRS_A|)
	\nonumber\\
	&\le & \sum_{\substack{A\in \SCRE,~c(A)=s\\ |A| \text{ is odd}}}m^s
	\nonumber\\
	&\le & \phi_sm^s.
}
\claimend

Now, by (\ref{epr-ii}) and Claims~\ref{cl1-4-9},~\ref{cl1-4-10} and~\ref{cl1-4-11}, we have
\eqn{in2-17}
{
	P_{DP}(G,\cov)-P(G,m)&=&\sum_{s=1}^{n-r_0+1}\sum_{\substack{A\in \SCRE\\ c(A)=s}}(-1)^{|A|}(|\SCRS_A|-m^{s})
	\nonumber\\
	&\le &-m^{n-r_0+1}+\sum_{s=1}^{n-r_0}\phi_sm^s,
}
where the inequality holds when $m>k$. As $k,\phi_1,\cdots,\phi_{n-r_0}$ are independent of the value of $m$, there exists an $M\in\N$, such that $P_{DP}(G,\cov)-P(G,m)<0$ for all $m\ge M$. Hence the result is proven.
\proofend

We shall conclude this section by proving Corollary~\ref{coroPDP<}.

\noindent\textit{Proof of Corollary~\ref{coroPDP<}:}
Let $E^*=\{e_1,e_2,\cdots,e_k\}\subseteq E_G(V_1,V_2)$,
where $e_i=u_iv_i$ 
with $u_i\in V_1$ and $v_i\in V_2$ 
for all $i\in [k]$.
For each $i\in [k]$, 
let $\overrightarrow{e_i}$ be the directed edge $(u_i,v_i)$,
and let $\overrightarrow{E^*}=
\{\overrightarrow{e_i}: i\in [k]\}$.
By Theorem~\ref{theoPDP<}, 
it suffices to verify that 
$\overrightarrow{E^*}$ is balanced 
on every cycle $C$ of $G$ 
with $|E(C)|<r_0$.

Let $C$ be any cycle of $G$ such that 
$|E(C)|<r_0$ and $|E(C)\cap E^*|$ is positive.
By the definition of $r_0$,
$|E(C)\cap E^*|=2r$ for some positive 
integer $r$, where $2r\le k$.
Without loss of generality, 
assume that 
$E(C)\cap E^*=\{e_i: i\in [2r]\}$.

Let $P$ be any minimal path of $C$ 
which contains exactly two edges in 
$E(C)\cap E^*$, say $e_i$ and $e_j$.
Obviously, by the minimality of $P$,
$e_i$ and $e_j$ must be the two edges
incident with two end-vertices of $P$. 
Then, the consecutive vertices on $P$ 
cannot appear 
in any one of the following 
orders:
$$
u_i,v_i,\cdots, u_j,v_j
\qquad 
\mbox{or}\qquad
v_i,u_i,\cdots, v_j,u_j.
$$
Otherwise, some component of $C-(E(C)\cap E^*)$ is either a $(v_i,u_j)$-path or a $(u_i,v_j)$-path 
in $G-E^*$,
contradicting the given condition
in the corollary.
Thus, the consecutive vertices on $P$ 
must appear in  one of the following 
orders:
$$
u_i,v_i,\cdots, v_j,u_j
\qquad 
\mbox{or}\qquad
v_i,u_i,\cdots, u_j,v_j.
$$
Since $|E(C)\cap E^*|=2r$, 
by the  definition of 
directed edges in  $\overrightarrow{E^*}$,
the above conclusion implies that 
the directed edges of $\overrightarrow{E^*}$ are balanced on $C$.

The corollary then follows from 
Theorem~\ref{theoPDP<}.
\proofend

\section{Study on plane graphs\label{secplane}}

By Corollary~\ref{corogd4},
every plane near-triangulation is DP-good and thus belongs to $DP^*$.
In the following, we consider 
those plane graphs $G$ in which 
at least two faces are not bounded by 
$3$-cycles. 
We will first show that 
such a plane graph $G$ may belong to
$DP_<$ if some face of $G$ 
is bounded by a $4$-cycle.

\begin{coro}\label{coro3PDP<}
	Let $G$ be any $2$-connected plane graph in which each $3$-cycle is the boundary of some face of $G$.
	If at least two faces of $G$ are not bounded by $3$-cycles 
	and one of them 
	is bounded by a $4$-cycle,  then $G\in DP_<$.
\end{coro}

\begin{figure}[!ht]
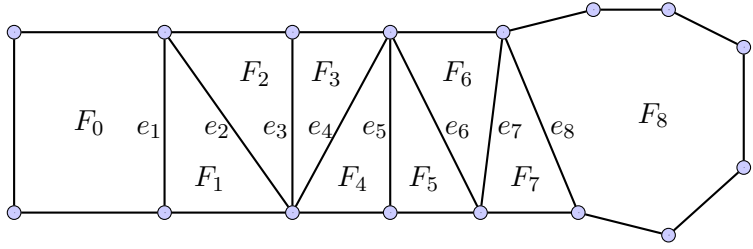

	\tikzstyle{cblue}=[circle, draw, thin,fill=blue!20, scale=0.5]
\tikzp{1}
{
	\foreach \place/\y in {{(0.5,1.2)/1},{(0.5,-1.2)/2},{(2.5,1.2)/3},{(2.5,-1.2)/4},{(4.2,1.2)/5},{(4.2,-1.2)/6},{(5.5,1.2)/7},{(5.5,-1.2)/8},{(7,1.2)/9},{(6.7,-1.2)/10},{(8,-1.2)/11},{(8.2,1.5)/12},{(9.2,1.5)/13},{(10.2,1)/14},{(10.2,-0.6)/15},{(9.2,-1.5)/16}}   
	\node[cblue] (b\y) at \place {};

	\filldraw[black] (b1) circle (0pt)node[anchor=south] {};
	\filldraw[black] (b2) circle (0pt)node[anchor=south] {};
	\filldraw[black] (b3) circle (0pt)node[anchor=south] {};
	\filldraw[black] (b4) circle (0pt)node[anchor=north] {};		
	\filldraw[black] (b5) circle (0pt)node[anchor=north] {};		
	\filldraw[black] (b6) circle (0pt)node[anchor=north] {};				
	\filldraw[black] (b7) circle (0pt)node[anchor=south] {};
	\filldraw[black] (b8) circle (0pt)node[anchor=south] {};
	\filldraw[black] (b9) circle (0pt)node[anchor=south] {};
	\filldraw[black] (b10) circle (0pt)node[anchor=north] {};		
	\filldraw[black] (b11) circle (0pt)node[anchor=north] {};		
	\filldraw[black] (b12) circle (0pt)node[anchor=north] {};				
	\filldraw[black] (b13) circle (0pt)node[anchor=south] {};
	\filldraw[black] (b14) circle (0pt)node[anchor=north] {};		
	\filldraw[black] (b15) circle (0pt)node[anchor=north] {};		
	\filldraw[black] (b16) circle (0pt)node[anchor=north] {};			
	
	\draw[thick]  (b1) -- (b2) -- (b4) -- (b3) -- (b1);	
	\draw[thick]  (b4) -- (b6) -- (b3) -- (b5) -- (b6);
	\draw[thick]  (b5) -- (b7) -- (b6) -- (b8) -- (b7);
	\draw[thick]  (b8) -- (b10) -- (b7) -- (b9) -- (b10);
	\draw[thick]  (b10) -- (b11) -- (b9) -- (b12) -- (b13) -- (b14) -- (b15) -- (b16) -- (b11);
	
	\node [style=none] (cap1) at (1.5,0)  {$F_0$};
	\node [style=none] (cap1) at (2.3,-0.1)  {$e_1$};
	\node [style=none] (cap1) at (3.1,-0.75)  {$F_1$};	
	\node [style=none] (cap1) at (3.2,-0.1)  {$e_2$};
	\node [style=none] (cap1) at (3.7,0.65)  {$F_2$};	
	\node [style=none] (cap1) at (3.98,-0.1)  {$e_3$};
	\node [style=none] (cap1) at (4.65,0.65)  {$F_3$};	
	\node [style=none] (cap1) at (4.58,-0.1)  {$e_4$};
	\node [style=none] (cap1) at (5,-0.75)  {$F_4$};	
	\node [style=none] (cap1) at (5.3,-0.1)  {$e_5$};
	\node [style=none] (cap1) at (5.95,-0.75)  {$F_5$};	
	\node [style=none] (cap1) at (6.4,-0.1)  {$e_6$};
	\node [style=none] (cap1) at (6.4,0.65)  {$F_6$};	
	\node [style=none] (cap1) at (7.1,-0.1)  {$e_7$};
	\node [style=none] (cap1) at (7.3,-0.75)  {$F_7$};	
	\node [style=none] (cap1) at (7.8,-0.1)  {$e_8$};
	\node [style=none] (cap1) at (9,0.1)  {$F_8$};	
}
	\caption{The graph consisting of vertices and edges on the boundaries of faces
		$F_0$, $F_1,\cdots,F_8$
}
	\label{remark-f1}
\end{figure}

\proof We can choose 
a shortest sequence of faces
$F_0,F_1,\cdots,F_t$ in $G$,
where $t\ge 1$, 
$F_0$ is bounded by a $4$-cycle
and $F_t$ is bounded by more than $3$ edges, 
such that $F_i$ is bounded by a $3$-cycle for each $i\in [t-1]$, and 
faces $F_{i-1}$ and $F_i$ share 
an edge $e_i$ on their boundaries
for each $i\in [t]$.
An example of the subgraph consisting of vertices and edges on boundaries of 
faces $F_0,F_1,\cdots,F_t$  
is shown in Figure~\ref{remark-f1}, where $t=8$.

If $t=1$, then $\ell_G(e_1)=4$
and thus $G\in DP_<$ by Theorem~\ref{Dong22DP<}.
Now assume that $t\ge 2$.
As $F_1,F_2,\cdots, F_{t-1}$ 
are all bounded by $3$-cycles,
$e_{i}$ and $e_{i+1}$ have a common end-vertex for each $i\in [t-1]$.
Thus, $e_i$ can be written as $e_i=u_iv_i$ for all $i\in [t]$  
such that either $u_i=u_{i+1}$
(i.e., $u_i$ and $u_{i+1}$ 
are the same vertex)
or $v_i=v_{i+1}$ for all $i\in [t-1]$.
Let $V_1=\{u_i:i\in [t]\}$ and $V_2=\{v_i:i\in [t]\}$. Then $E^*:=\{e_i:i\in [t]\}\subseteq E_G(V_1,V_2)$. 

As $F_0$ is bounded by a $4$-cycle, 
$G$ has a $4$-cycle $C$ with 
$|E(C)\cap E^*|=1$.
But, as each $3$-cycle in $G$ must be the boundary of some face of $G$,
there is no $3$-cycle $C$ in $G$ 
with $|E(C)\cap E^*|=1$.
As the dual edges of the edges in $E^*$ 
actually form a shortest path
connecting vertices $F_0^*$ and $F_t^*$
in the dual plane graph $G^*$ of $G$,
there is no $3$-cycle $C$ in $G$ 
with $|E(C)\cap E^*|=3$.
Therefore, $\g_G(E^*)=4$.
Thus, by Corollary~\ref{coro2PDP<}, 
$G\in DP_<$, and the result holds.
\proofend

It is not difficult to generalize Corollary~\ref{coro3PDP<}
as stated below.

\begin{coro}\label{coro4PDP<}
	Let $G$ be any $2$-connected plane graph. If 
	$F_0,F_1,\cdots, F_t$ are faces in $G$,
	where $t\ge 1$, which satisfy 
	the following conditions, 
	then  $G\in DP_<$:
	\begin{enumerate}
		\item  
		only $F_0$ and $F_t$ are not bounded by $3$-cycles,
		$F_0$ is bounded by an even cycle $C_r$ and $F_t$  is bounded by a cycle not shorter than $r$;
		\item for each $i\in [t]$,
		 faces $F_{i-1}$ and $F_i$ 
		 share an edge $e_i$ on their 
		  boundaries; and 
		\item for $E^*=\{e_i: i\in [t]\}$, if $C$ is a cycle in $G$ with $E(C)\cap E^*\ne \emptyset$, then 
		$|E(C)|\ge r$ holds whenever 
		either $F_0$ or $F_t$ is within
		cycle $C$.
	\end{enumerate}
\end{coro}

\begin{figure}[!ht]
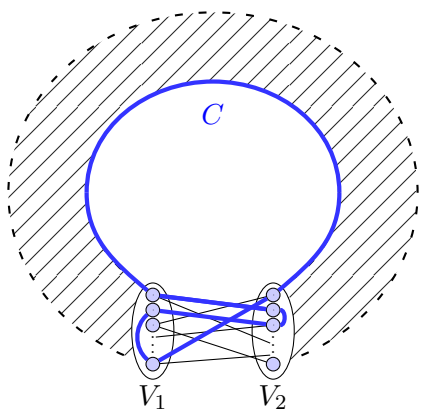

	\tikzstyle{cblue}=[circle, draw, thin,fill=blue!20, scale=0.5]
	\tikzp{0.8}
	{
		\draw[thick, pattern= my north east lines, pattern color=black!85, loosely dashed](0,0) ellipse (3.4 and 3);
		\draw[thick, loosely dashed, fill=white](0,0) ellipse (2.1 and 1.85);
		\draw[thick, white, fill=white](0,-2.5) ellipse (1.2 and 1.2);		
		\draw[ fill=white] (-1,-2.3) ellipse (0.35 and 0.8);
		\draw[ fill=white] (1,-2.3) ellipse (0.35 and 0.8);

		\foreach \place/\y in {{(-1,-1.7)/1},{(1,-1.7)/6},{(-1,-1.95)/3},{(1,-1.95)/2},{(-1,-2.2)/5},{(1,-2.2)/4},{(-1,-2.85)/7},{(1,-2.85)/8}}   
		\node[cblue] (b\y) at \place {};

		\filldraw[black] (b1) circle (0pt)node[anchor=south] {};
		\filldraw[black] (b6) circle (0pt)node[anchor=south] {};
		\filldraw[black] (b3) circle (0pt)node[anchor=south] {};
		\filldraw[black] (b2) circle (0pt)node[anchor=north] {};		
		\filldraw[black] (b5) circle (0pt)node[anchor=north] {};		
		\filldraw[black] (b4) circle (0pt)node[anchor=north] {};

	\draw[] (0.95,-2.5) -- (b1); 
\draw[]  (b5) -- (b6);	
\draw[] (-0.95,-2.4) -- (b4) ;
\draw[ ] (b3) -- (b4);	 
\draw[ ] (b5) -- (b8);
\draw[ ] (b7) -- (0.95,-2.7);
\draw[ blue!80, ultra thick] 	(b2)  to [bend left=90]  (b4) ;
\draw[ blue!80, ultra thick] 	(b7)  to [bend left=50]  (b3) ;
\draw[ blue!80, ultra thick] 	(b1) -- (b2);
\draw[ blue!80, ultra thick] 	(b3) -- (b4);
\draw[ blue!80, ultra thick] 	(b1) -- (b2);
\draw[ blue!80, ultra thick] 	(b7) -- (b6);
\draw[blue!80, ultra thick] (b1) to [out=140,in=270]  (-2.1,0);
\draw[blue!80, ultra thick] (-2.1,0) to [out=90,in=180]  (0,1.85);
\draw[blue!80, ultra thick] (b6) to [out=40,in=270]  (2.1,0);
\draw[blue!80, ultra thick] (2.1,0) to [out=90,in=0]  (0,1.85);

		\node [style=none] (cap1) at (0,1.3)  {$\blue{C}$};
		\node [style=none] (cap2) at (-1,-3.4)  {$V_1$};
		\node [style=none] (cap3) at (1,-3.4)  {$V_2$};
		\node [style=none] (cap4) at (-1,-2.43)  {\tiny{$\vdots$}};
		\node [style=none] (cap5) at (1,-2.43)  {\tiny{$\vdots$}};
	}
	
	\caption{
		$E^*\subseteq E_G(V_1,V_2)$ for 
	 $V_1,V_2\subseteq V(G)$ and $C\in \C'_G(E^*)$}
	\label{fig1-5}
\end{figure}

By Corollarys~\ref{coroPDP<} and~\ref{coro3PDP<},
it is interesting to notice that 
quite many graphs in $DP_<$ have 
their structures with 
a doughnut shape,
as shown in Figure~\ref{fig1-5},
where $E^*\subseteq E_G(V_1,V_2)$ for two  disjoint vertex sets $V_1$ and $V_2$,
and $C$ is a shortest cycle in $G$
such that $|E(C)\cap E^*|$ is odd
and $|E(C)|$ is even.

However, for some other plane graphs
which also look like doughnuts, we still don't know  whether they belong to $DP_{\approx}$ or $DP_{<}$.
For example, for a $2$-connected plane graph $G$ which is not a near-triangulation,  
if $\ell_G(e)=3$ for all $e\in E(G)$,
and those faces in $G$ not bounded by 
$3$-cycles have respectively  $q_1,q_2,\cdots,q_t$ edges on 
their boundaries, 
where $4\le q_1\le q_2\le \cdots \le q_t$ and $q_j$ is even 
whenever $q_i<q_j$ and $q_i$ is even,
it is  still 
unknown if $G$ belongs to  $DP_{\approx}$ or 
$DP_{<}$.
For the particular case that 
$q_i$ is odd for all $i\in [t-1]$, 
we guess $G$ belongs to $DP_{\approx}$.

\end{document}